\pdfoutput=1 

\documentclass[11pt]{amsart}

\usepackage{amssymb,amsthm,enumitem,colonequals,tikz-cd,microtype,quiver}
\usepackage[normalem]{ulem} 
\usepackage{mathtools} 
\usepackage{graphicx}

\usepackage[osf]{Baskervaldx}
\usepackage[baskervaldx]{newtxmath}
\usepackage[cal=boondoxo]{mathalfa} 

\usepackage[top=3.75cm, bottom=3cm, left=3.5cm, right=3.5cm]{geometry}

\usepackage{xcolor}
\colorlet{darkblue}{blue!55!black}
\colorlet{darkcyan}{cyan!50!black}
\colorlet{darkgreen}{green!60!black}

\PassOptionsToPackage{hyphens}{url}
\usepackage{hyperref}
\hypersetup{
colorlinks=true, 
linkcolor=darkblue,
urlcolor=darkcyan,
citecolor=darkgreen,
}

\def\eqref#1{\textcolor{darkblue}{(\ref{#1})}}

\usepackage[nameinlink,capitalize]{cleveref} 
\Crefformat{section}{#2\S#1#3}
\Crefmultiformat{section}{#2\S\S#1#3}{ and~#2#1#3}{, #2#1#3}{, and~#2#1#3}

\usepackage[pagewise]{lineno}
\let\oldequation\equation
\let\oldendequation\endequation
\renewenvironment{equation}{\linenomathNonumbers\oldequation}{\oldendequation\endlinenomath}
\expandafter\let\expandafter\oldequationstar\csname equation*\endcsname
\expandafter\let\expandafter\oldendequationstar\csname endequation*\endcsname
\renewenvironment{equation*}{\linenomathNonumbers\oldequationstar}{\oldendequationstar\endlinenomath}
\let\oldalign\align
\let\oldendalign\endalign
\renewenvironment{align}{\linenomathNonumbers\oldalign}{\oldendalign\endlinenomath}
\expandafter\let\expandafter\oldalignstar\csname align*\endcsname
\expandafter\let\expandafter\oldendalignstar\csname endalign*\endcsname
\renewenvironment{align*}{\linenomathNonumbers\oldalignstar}{\oldendalignstar\endlinenomath}

\theoremstyle{plain}
\newtheorem{theorem}{Theorem}[section]
\newtheorem{lemma}[theorem]{Lemma}
\newtheorem{corollary}[theorem]{Corollary}
\newtheorem{proposition}[theorem]{Proposition}

\theoremstyle{definition}
\newtheorem{definition}[theorem]{Definition}
\newtheorem{example}[theorem]{Example}
\newtheorem{remark}[theorem]{Remark}

\newtheorem{question}[theorem]{Question}
\newtheorem{remark-definition}[theorem]{Remark/Definition}

\newtheorem*{ack}{Acknowledgments}

\AddToHook{env/conjecture/begin}{\crefalias{theorem}{conjecture}}
\AddToHook{env/lemma/begin}{\crefalias{theorem}{lemma}}
\AddToHook{env/corollary/begin}{\crefalias{theorem}{corollary}}
\AddToHook{env/proposition/begin}{\crefalias{theorem}{proposition}}
\AddToHook{env/definition/begin}{\crefalias{theorem}{definition}}
\AddToHook{env/remark/begin}{\crefalias{theorem}{remark}}
\AddToHook{env/setup/begin}{\crefalias{theorem}{setup}}
\AddToHook{env/example/begin}{\crefalias{theorem}{example}}
\AddToHook{env/conjecture/begin}{\crefalias{theorem}{conjecture}}
\AddToHook{env/notation/begin}{\crefalias{theorem}{notation}}
\AddToHook{env/reminder/begin}{\crefalias{theorem}{reminder}}
\AddToHook{env/question/begin}{\crefalias{theorem}{question}}
\AddToHook{env/remark-definition/begin}{\crefalias{theorem}{remark-definition}}

\makeatletter
\def\slashedarrowfill@#1#2#3#4#5{%
$\m@th\thickmuskip0mu\medmuskip\thickmuskip\thinmuskip\thickmuskip
\relax#5#1\mkern-7mu%
\cleaders\hbox{$#5\mkern-2mu#2\mkern-2mu$}\hfill
\mathclap{#3}\mathclap{#2}%
\cleaders\hbox{$#5\mkern-2mu#2\mkern-2mu$}\hfill
\mkern-7mu#4$%
}
\def\rightslashedarrowfill@{%
\slashedarrowfill@\relbar\relbar\mapstochar\rightarrow}
\newcommand\xslashedrightarrow[2][]{%
\ext@arrow 0055{\rightslashedarrowfill@}{#1}{#2}}
\makeatother

\numberwithin{equation}{section}
\numberwithin{theorem}{section}

\title[Krull dimension is a Fourier--Mukai invariant]{Krull dimension is a Fourier--Mukai invariant}

\author[P.~Lank]{Pat Lank}
\address{P.~Lank,
Dipartimento di Matematica “F. Enriques”, Universit\`{a} degli Studi di Milano, Via Cesare
Saldini 50, 20133 Milano, Italy}
\email{lankp.math@gmail.com}

\date{\today}



\begin{document}

\begin{abstract}
	We prove that Krull dimension is invariant under Fourier--Mukai partnerships between algebraic spaces which are proper over an arbitrary field.
	Moreover, if one of the partners is Cohen--Macaulay, then so is the other.
	In addition, we prove the relative version of these results.
\end{abstract}

\maketitle

\tableofcontents

\section{Introduction}
\label{sec:intro}

\subsection{Motivation}
\label{sec:intro_motivation}

The bounded derived category of coherent sheaves encodes information about the
geometry of an algebraic space.
This is particularly prominent in the study of singularities.
It is natural to ask which local properties can be recovered from this category.

We consider the following question.

\begin{question}
	\label{q:invariants}
	Which invariants of the local rings of a proper algebraic space over a field
	are preserved under Fourier--Mukai equivalences?
\end{question}

Let $k$ be a field.
Suppose $Y_1$ and $Y_2$ are proper algebraic spaces over $k$.
For $E\in D^b_{\operatorname{coh}}(Y_1\times_kY_2)$, write $\Phi_E$ for the
integral transform with kernel $E$.
See \Cref{sec:prelim_integral_transforms}.
We call $Y_1$ and $Y_2$ \textit{Fourier--Mukai partners} if, for such an
$E$, the functor $\Phi_E$ restricts to an equivalence
\begin{displaymath}
	\Phi_E\colon D^b_{\operatorname{coh}}(Y_1)
	\to D^b_{\operatorname{coh}}(Y_2).
\end{displaymath}
\Cref{q:invariants} is about invariants preserved by equivalences arising from
integral transforms, rather than arbitrary equivalences of triangulated
categories.

For projective $k$-schemes, every exact $k$-linear equivalence is naturally isomorphic to a functor of the form $\Phi_E$ \cite[Corollary 1.3]{Ballard:2009}.
Hence, for projective $k$-schemes, our results apply to arbitrary exact $k$-linear equivalences.
See \Cref{cor:intro_projective}. 
Our main result concerns Krull dimension and requires none of smoothness, projectivity, or equidimensionality.

\subsection{Background}
\label{sec:intro_background}

For a smooth projective scheme $Y$ of pure Krull dimension $d$ over a field, the Serre functor on $D^b_{\operatorname{coh}}(Y)$ is given by tensoring with the canonical bundle and shifting by $d$.
An equivalence intertwines Serre functors, providing a familiar explanation
for the invariance of Krull dimension in the smooth projective setting.
See e.g.\ \cite[Theorem 2.3(0)]{Kawamata:2002}.

In the presence of singularities, the dualizing complex replaces the shifted
canonical bundle.
However, it is locally isomorphic to a shift of an invertible sheaf only around Gorenstein points, where the shift may vary between connected components. 
In particular, the same argument does not directly recover the Krull dimension of the algebraic space, and so different strategies are required.

Hern\'{a}ndez Ruip\'{e}rez, L\'{o}pez Mart\'{\i}n, and Sancho de Salas first established preservation results in the singular setting.
For projective schemes over an algebraically closed field, \cite[Theorem 4.4]{Ruiperez/Hernandez/Martin/SanchodeSalas:2009} shows that a reduced Fourier--Mukai partner of an equidimensional Cohen--Macaulay scheme is equidimensional of the same dimension, that an equidimensional partner of the same dimension is Cohen--Macaulay, and that such a partner is Gorenstein whenever the given scheme is.
For Gorensteinness, \cite[Theorem 1.3]{GuisadoVillaalgordo/Lank:2026b} removes the hypotheses of projectivity and equidimensionality.

The present manuscript considers proper algebraic spaces over an arbitrary
field, which may have arbitrary singularities. 
Our focus is on the depths at closed points.
The argument of \cite{Ruiperez/Hernandez/Martin/SanchodeSalas:2009} uses local complete intersection zero-cycles.
Instead, we use perfect complexes supported at closed points, and track only the top nonzero cohomological degree of the dualizing complex.
In particular, \Cref{prop:cm_invariance} recovers the Cohen--Macaulay assertions of \cite[Theorem 4.4]{Ruiperez/Hernandez/Martin/SanchodeSalas:2009} without assuming projectivity, an algebraically closed base field, reducedness, or equidimensionality of the partner.

\subsection{Results}
\label{sec:intro_results}

Our main theorem is as follows.

\begin{theorem}
	\label{introthm:main}
	Let $Y_1$ and $Y_2$ be proper algebraic spaces over a field $k$ which are
	Fourier--Mukai partners.
	Then
	\begin{displaymath}
		\dim Y_1=\dim Y_2.
	\end{displaymath}
	Moreover, there exists a bijection between the connected components of $|Y_1|$ and $|Y_2|$ such that corresponding components are Fourier--Mukai partners, have the same Krull dimension, and have the same set of depths at closed points.
	Furthermore, $Y_1$ is Cohen--Macaulay if, and only if, $Y_2$ is Cohen--Macaulay.
	In fact, the same holds for being Cohen--Macaulay and equidimensional.
\end{theorem}

\begin{proof}
	This follows from \Cref{cor:dimension_invariance,prop:connected_components,prop:cm_invariance}.
\end{proof}

\begin{corollary}
	\label{cor:intro_projective}
	Let $Y_1$ and $Y_2$ be projective schemes over a field $k$.
	Suppose there exists an exact $k$-linear equivalence $D^b_{\operatorname{coh}}(Y_1)\to D^b_{\operatorname{coh}}(Y_2)$.
	Then the conclusions of \Cref{introthm:main} hold for $Y_1$ and $Y_2$.
	In particular, $\dim Y_1=\dim Y_2$.
\end{corollary}

\begin{proof}
	By \cite[Corollary 1.3]{Ballard:2009}, the equivalence is naturally isomorphic to the restriction of $\Phi_E$ for some $E\in D^b_{\operatorname{coh}}(Y_1\times_k Y_2)$.
	Hence, $Y_1$ and $Y_2$ are Fourier--Mukai partners, and so \Cref{introthm:main} applies.
\end{proof}

For instance, the comparison of mixed Hodge numbers of singular complex projective Fourier--Mukai partners in \cite[Theorem 1.1]{Konovalov:2026} is stated for partners of the same dimension.
Now, by \Cref{cor:intro_projective}, this is automatic.
However, this requires working in the enhanced setting, whereas Fourier--Mukai partnerships are defined here in the triangulated setting.
See \cite[Appendix C]{GuisadoVillaalgordo/Lank:2026} for limitations.

Moreover, there exist relative consequences.

\begin{corollary}
	\label{cor:intro_relative}
	Let $S$ be a Noetherian algebraic space and $d\geq0$.
	Suppose $f_i\colon Y_i\to S$ are proper flat morphisms of algebraic spaces.
	Assume that $Y_1$ and $Y_2$ are Fourier--Mukai $S$-partners.
	For every $s\in|S|$, write
	$(Y_i)_s=Y_i\times_S\operatorname{Spec}(\kappa(s))$.
	\begin{enumerate}
		\item \label{cor:intro_relative1} $\dim (Y_1)_s=\dim (Y_2)_s$ for every $s\in|S|$
		\item\label{cor:intro_relative2} $(Y_1)_s$ is Cohen--Macaulay if, and only if, $(Y_2)_s$ is Cohen--Macaulay; the same holds for being Cohen--Macaulay and
		equidimensional
		\item \label{cor:intro_relative3} $f_1$ has relative dimension $\leq d$ if, and only if, $f_2$ has relative dimension $\leq d$ (see \cite[\href{https://stacks.math.columbia.edu/tag/04NH}{Tag 04NH}]{StacksProject}).
		\item \label{cor:intro_relative4} $\dim Y_1=\dim Y_2$
		\item \label{cor:intro_relative5} $Y_1$ is Cohen--Macaulay if, and only if, $Y_2$ is Cohen--Macaulay.
	\end{enumerate}
\end{corollary}

\subsection{Comparisons}
\label{sec:intro_comparisons}

These results appear to be new even for proper schemes over a field with arbitrary singularities.
To put it in context, we compare to the literature:
\begin{itemize}
	\item For smooth proper varieties over a field, \cite[Theorem 2]{Olander:2021} proves the statement that a fully faithful exact $k$-linear functor $D^b_{\operatorname{coh}}(Y_1)\to D^b_{\operatorname{coh}}(Y_2)$ forces $\dim Y_1\leq\dim Y_2$.
	The analogous statement for smooth proper tame algebraic stacks is \cite[Theorem 1.3]{Peng:2024}.
	Proofs of \cite{Olander:2021,Peng:2024} differ from the present approach.
	\item For proper complex varieties with rational singularities, the inequality $\dim Y_1\leq\dim Y_2$ is proven in \cite[Proposition 4.4]{Dutta/Lank/ManaliRahul:2025} for fully faithful integral transforms which preserve perfect complexes.
	In fact, this holds more generally, see loc.\ cit.
	\item Lower bounds also come from Rouquier dimension: if $Y$ is a reduced separated scheme of finite type over $k$, then $\dim Y$ is at most the Rouquier dimension of $D^b_{\operatorname{coh}}(Y)$ \cite[Proposition 7.17]{Rouquier:2008}.
	As Rouquier dimension is invariant under exact equivalences, this bounds the dimension of a reduced Fourier--Mukai partner, but it does not give equality.
	Indeed, the Rouquier dimension of $D^b_{\operatorname{coh}}(Y)$ is not known to agree with $\dim Y$ in general.
	\item In special classes, reconstruction results show that a partner stays in the class, and often that it is isomorphic to the given space.
	For example, see \cite{Ballard:2011} for Gorenstein projective varieties with ample or anti-ample canonical bundle, \cite{Spence:2023} for integral projective curves over an algebraically closed field of characteristic zero, \cite{LopezMartin:2014} for reduced projective Gorenstein curves of genus one with trivial dualizing sheaf, and \cite{Matsukawa:2025} for proper partners of abelian varieties.
\end{itemize}

None of these results applies to a proper algebraic space with arbitrary singularities, nor recover the equality of the sets of depths at closed points established below.

\subsection{Methodology}
\label{sec:intro_methods}

\subsubsection{Key step}
\label{sec:intro_methods_key_step}

The essential ingredient is stronger than invariance of Krull dimension. 
For each closed point $x\in|Y_1|$, there exists a nonzero perfect complex
$P_x$ on $Y_1$ with $\operatorname{Supp}(P_x)=\{x\}$.
If $\Phi_E$ induces the equivalence, then \Cref{thm:depth_correspondence}
shows that every closed point $y\in\operatorname{Supp}(\Phi_E(P_x))$ satisfies
\begin{equation}
	\label{eq:intro_depth_relation}
	\operatorname{depth}_x(Y_1)=\operatorname{depth}_y(Y_2).
\end{equation}
Furthermore, each $\Phi_E(P_x)$ has a closed point in its support, and every
closed point of $Y_2$ belongs to the support of some $\Phi_E(P_x)$.
Consequently,
\begin{equation}
	\label{eq:intro_depth_spectra}
	\{\operatorname{depth}_x(Y_1):x\in|Y_1|\textrm{ closed}\}
	=\{\operatorname{depth}_y(Y_2):y\in|Y_2|\textrm{ closed}\}.
\end{equation}
Thus, the set of depths of closed points is a Fourier--Mukai
invariant.
See \Cref{thm:depth_correspondence}.
A crucial step is that for every nonempty algebraic space $X$ finitely presented over a field,
\begin{displaymath}
	\dim X=\max\{\operatorname{depth}_x(X):x\in|X|\textrm{ closed}\}.
\end{displaymath}
See \Cref{lem:max_closed_depth}.
In particular, no equidimensionality hypothesis is required. 

The assertion for being Cohen--Macaulay and equidimensional follows from a
characterization in terms of depths at closed points.
We establish this in the proof of \Cref{prop:cm_invariance}.
To remove equidimensionality, we show that Fourier--Mukai partnerships match connected components.
See \Cref{prop:connected_components}.
This uses orthogonal decompositions of $D^b_{\operatorname{coh}}$. 
Unlike semiorthogonal decompositions, such decompositions into indecomposable pieces are unique up to reordering, see \Cref{lem:orthogonal_unique}.
A key fact is that connected Cohen--Macaulay algebraic spaces of finite presentation over a field are equidimensional.
See \Cref{lem:connected_cm_equidimensional}.

\subsubsection{Argument overview}
\label{sec:intro_methods_argument_overview}

We highlight how \eqref{eq:intro_depth_spectra} is established.
Write $f_i\colon Y_i\to\operatorname{Spec}(k)$ for the structure morphisms.
The proof of \Cref{thm:depth_correspondence} compares the top cohomological
degrees of two descriptions of $\Phi_E(P_x\otimes^{\mathbf{L}} f_1^!\mathcal{O}_{\operatorname{Spec}(k)})$
at a closed point $y\in\operatorname{Supp}(\Phi_E(P_x))$.

The first description comes from Serre duality.
For a proper algebraic space $Y$ over $k$, Serre duality describes
$\operatorname{Hom}_Y(Q,P\otimes^{\mathbf{L}}
	f^!\mathcal{O}_{\operatorname{Spec}(k)})$ for $P$ perfect and
$Q\in D^b_{\operatorname{coh}}(Y)$ in terms of the $k$-linear dual of
$\operatorname{Hom}_Y(P,Q)$.
See \Cref{lem:serre_duality,lem:transport_duality}.
By \cite[Remark 2.7 and Corollary 1.4]{GuisadoVillaalgordo/Lank:2026}, integral transforms inducing Fourier--Mukai partnerships preserve perfect complexes.
Hence, the associated duality isomorphisms yield
\begin{equation}
	\label{eq:intro_duality_transport}
	\Phi_E(P_x\otimes^{\mathbf{L}}
	f_1^!\mathcal{O}_{\operatorname{Spec}(k)})
	\cong\Phi_E(P_x)\otimes^{\mathbf{L}}
	f_2^!\mathcal{O}_{\operatorname{Spec}(k)}.
\end{equation}
See \Cref{lem:transport_duality}.

For the second description, \Cref{lem:point_tensor_comparison} gives a
local calculation of the top cohomological degree:
\begin{displaymath}
	\begin{aligned}
		&\operatorname{sup}_y (\Phi_E(P_x\otimes^{\mathbf{L}}
			f_1^!\mathcal{O}_{\operatorname{Spec}(k)}) )\\
		&\qquad=\operatorname{sup}_y(\Phi_E(P_x))
			+\operatorname{sup}_x(f_1^!\mathcal{O}_{\operatorname{Spec}(k)}).
	\end{aligned}
\end{displaymath}
The relevant pushforward calculation is possible because the complexes
in question are supported over the closed point $x$, whose residue field
is finite over $k$.
More precisely, \Cref{lem:finite_support_exact} proves that pushforward
along a proper morphism commutes with cohomology for bounded pseudocoherent
complexes supported on a closed subspace finite over the base, and
preserves their top nonzero cohomological degree.

At a closed point $z$ of $Y_i$, the complex
$f_i^!\mathcal{O}_{\operatorname{Spec}(k)}$ is locally a normalized
dualizing complex.
Since a normalized dualizing complex has vanishing cohomology in positive degrees, its top nonzero cohomological degree is the negative of the depth.
See \cite[\href{https://stacks.math.columbia.edu/tag/0A7U}{Tags 0A7U} and \href{https://stacks.math.columbia.edu/tag/0AWR}{0AWR}]{StacksProject}.
It follows that $\operatorname{sup}_z(f_i^!\mathcal{O}_{\operatorname{Spec}(k)})=-\operatorname{depth}_z(Y_i)$. 
On the right-hand side of \eqref{eq:intro_duality_transport},
\Cref{lem:pointwise_tensor} expresses the top degree as
$\operatorname{sup}_y(\Phi_E(P_x))-\operatorname{depth}_y(Y_2)$.
Comparing both sides therefore gives \eqref{eq:intro_depth_relation}.
Finally, a `detection property' of point supported perfect complexes
(see \Cref{lem:point_perfect}) shows that the supports of their transforms
account for every closed point of $Y_2$.
This yields \eqref{eq:intro_depth_spectra}.

\begin{ack}
	ERC Advanced Grant 101095900-TriCatApp for support.
	Large language models were used to assist with proofreading, refining the exposition, checking arguments, and exploring proof strategies.
\end{ack}

\section{Preliminaries}
\label{sec:preliminaries}

We follow \cite{StacksProject} for conventions on algebraic spaces.
We call a morphism \textbf{natural} if it is functorial, and \textbf{canonical} if it is given by an explicit construction that need not be functorial.
Let $Y$ be a quasi-compact algebraic space.
An \textbf{\'{e}tale presentation} of $Y$ is a morphism $U \to Y$ from a scheme which is \'{e}tale, surjective, and finitely presented.
We follow \cite[\S 5]{Lank:2026} for the notion of \textbf{germs} on algebraic spaces.

\subsection{Notation}
\label{sec:prelim_notation}

We work with the small \'{e}tale ringed site $(Y_{\textrm{et}},\mathcal{O}_Y)$.
$D(\mathcal{O}_Y)$ is the unbounded derived category of sheaves of $\mathcal{O}_Y$-modules.
$D_{\operatorname{qc}}(Y)$ is the strictly full subcategory of $D(\mathcal{O}_Y)$ consisting of complexes with quasi-coherent cohomology.
For $?\in \{b,+,-\}$, let $D^{?}_{\operatorname{qc}}(Y)$ be the strictly full subcategory of $D_{\operatorname{qc}}(Y)$ with prescribed cohomological vanishing.
$\operatorname{Perf}(Y)$ is the strictly full subcategory of $D_{\operatorname{qc}}(Y)$ consisting of perfect complexes.
$\operatorname{Qcoh}(Y)$ is the abelian category of quasi-coherent sheaves.
If $Y$ is Noetherian, $\operatorname{coh}(Y)$ is the abelian category of coherent sheaves.
We use the notion of \textbf{pseudocoherent} complexes on the ringed site $(Y_{\textrm{et}},\mathcal{O}_Y)$ from \cite[\href{https://stacks.math.columbia.edu/tag/08FT}{Tag 08FT}]{StacksProject}.
If $Y$ is Noetherian, an object of $D_{\operatorname{qc}}(Y)$ is pseudocoherent if, and only if, it lies in $D^-_{\operatorname{coh}}(Y)$ \cite[\href{https://stacks.math.columbia.edu/tag/08IK}{Tag 08IK}]{StacksProject}, and so $D^b_{\operatorname{coh}}(Y)$ is the strictly full subcategory of $D_{\operatorname{qc}}(Y)$ consisting of the pseudocoherent complexes with bounded cohomology.

Set $H^i(Y_{\textrm{et}},-)$ for $H^i\mathbf{R}\Gamma(Y_{\textrm{et}},-)$.
Recall that the derived internal Hom $\operatorname{\mathbb{R}\mathcal{H}\! \mathit{om}}_{\mathcal{O}_Y}$ is formed in the ambient category $D(\mathcal{O}_Y)$ \cite[\href{https://stacks.math.columbia.edu/tag/08J7}{Tag 08J7}]{StacksProject}.
In particular, it need not agree with the right adjoint to tensor product in $D_{\operatorname{qc}}(Y)$, which is denoted by $\operatorname{\mathbf{R}\mathcal{H}\! \mathit{om}}$.
See \cite[Lemma 2.11]{GuisadoVillaalgordo/Lank:2026}.
Throughout, we only use the former if its first argument is perfect and its second argument is in $D_{\operatorname{qc}}(Y)$.
In this case, by \cite[\href{https://stacks.math.columbia.edu/tag/0A8A}{Tag 0A8A}]{StacksProject}, the complex lies in $D_{\operatorname{qc}}(Y)$.
For $K \in \operatorname{Perf}(Y)$, let $K^\vee:=\operatorname{\mathbb{R}\mathcal{H}\! \mathit{om}}_{\mathcal{O}_Y}(K,\mathcal{O}_Y)$.
By \cite[\href{https://stacks.math.columbia.edu/tag/08JJ}{Tag 08JJ}]{StacksProject}, $K^\vee$ is perfect, $(K^\vee)^\vee\cong K$, and for every $M\in D(\mathcal{O}_Y)$ there exist natural isomorphisms
\begin{equation}
	\label{eq:perfect_dual}
	M\otimes^{\mathbf{L}}_{\mathcal{O}_Y}K^\vee
	\cong\operatorname{\mathbb{R}\mathcal{H}\! \mathit{om}}_{\mathcal{O}_Y}(K,M),
	\qquad
	H^0(Y_{\textrm{et}},M\otimes^{\mathbf{L}}_{\mathcal{O}_Y}K^\vee)
	\cong\operatorname{Hom}_Y(K,M).
\end{equation}

\subsection{Krull dimension}
\label{sec:prelim_dimension}

Denote by $|Y|$ the topological space of $Y$ \cite[\href{https://stacks.math.columbia.edu/tag/03BT}{Tag 03BT}]{StacksProject}.

\begin{lemma}
	\label{lem:closed_points_charts}
	Let $X$ be a decent algebraic space (e.g.\ quasi-separated \cite[\href{https://stacks.math.columbia.edu/tag/03JX}{Tag 03JX}]{StacksProject}).
	Pick an \'{e}tale morphism $a\colon U\to X$ with $U$ an affine scheme.
	If $u\in U$ maps to a closed point of $|X|$, then $u$ is a closed point of $|U|$.
\end{lemma}

\begin{proof}
	Since $U$ is quasi-compact, the closure of $\{u\}$ in $U$ contains a point $u^\prime$ which is closed in $U$ \cite[\href{https://stacks.math.columbia.edu/tag/005E}{Tag 005E}]{StacksProject}.
	Then $a(u^\prime)$ lies in the closure of $\{a(u)\}$.
	Hence, by hypothesis, the closure of $\{a(u)\}$ is $\{a(u)\}$.
	Thus, $u$ and $u^\prime$ map to the same point of $|X|$, and $u$ specializes to $u^\prime$.
	Therefore, \cite[\href{https://stacks.math.columbia.edu/tag/03K5}{Tag 03K5}]{StacksProject} implies $u=u^\prime$.
\end{proof}

Let $X$ be an algebraic space.
Choose $x\in|X|$.
Recall that the \textbf{dimension of $X$ at $x$}, denoted $\dim_x(X)$, is the dimension $\dim_u(U)$ of $U$ at $u$ for any (equivalently, some) \'{e}tale morphism $a\colon U\to X$ from a scheme and any $u\in U$ such that $a(u)=x$ \cite[\href{https://stacks.math.columbia.edu/tag/04N5}{Tag 04N5}]{StacksProject}.
Here, $\dim_u (U)$ is the Krull dimension of $|U|$ at $u$; see \cite[\href{https://stacks.math.columbia.edu/tag/0055}{Tag 0055}]{StacksProject} and the proof of \cite[\href{https://stacks.math.columbia.edu/tag/04N4}{Tag 04N4}]{StacksProject}.
The \textbf{Krull dimension} of $X$ is given by
\begin{displaymath}
	\dim X=\sup_{x\in|X|}\dim_x(X).
\end{displaymath}
See \cite[\href{https://stacks.math.columbia.edu/tag/04N6}{Tag 04N6}]{StacksProject}.

If $X$ is a scheme, Krull dimension agrees with the supremum of Krull dimensions of its local rings at all points in $|X|$ \cite[\href{https://stacks.math.columbia.edu/tag/04MU}{Tag 04MU}]{StacksProject}.
If $X$ is decent, then $\dim_x(X)=\dim_x(|X|)$ and $\dim X=\dim|X|$ \cite[\href{https://stacks.math.columbia.edu/tag/0A4J}{Tag 0A4J}]{StacksProject}.
That is, Krull dimension defined in terms of \'{e}tale germs agrees with the topological analog.
In general, the dimension of $X$ at $x$ should not be confused with the Krull dimension of the local ring of $X$ at $x$ (see \Cref{def:local-depth-dimension}), e.g.\ see the `Warning' of \cite[\href{https://stacks.math.columbia.edu/tag/04N3}{Tag 04N3}]{StacksProject}.
By \Cref{lem:dimension_charts}, these agree at closed points if $X$ is finitely presented over a field.

\begin{definition}
	\label{def:equidimensional}
	Let $X$ be an algebraic space locally of finite type over a field $k$.
	We say that $X$ is \textbf{equidimensional} if $\dim_x(X)=\dim X$ for every $x\in|X|$.
\end{definition}

\begin{remark}
	\label{rmk:equidimensional_schemes}
	Let $X$ be a scheme locally of finite type over a field $k$.
	For every $x\in X$, $\dim_x(X)$ is the maximum of the Krull dimensions of the irreducible components of $X$ containing $x$ \cite[\href{https://stacks.math.columbia.edu/tag/0A21}{Tag 0A21}]{StacksProject}.
	Note that the generic point of an irreducible component lies in no other irreducible component.
	Hence, $X$ is equidimensional in the sense of \Cref{def:equidimensional} if, and only if, all irreducible components of $X$ have the same Krull dimension.
	Consequently, \Cref{def:equidimensional} agrees with the usual notion for schemes.
	Moreover, if $X$ is a decent algebraic space, then \cite[\href{https://stacks.math.columbia.edu/tag/0A4J}{Tag 0A4J}]{StacksProject} shows that the condition in \Cref{def:equidimensional} only depends on the topological space $|X|$.
\end{remark}

\begin{lemma}
	\label{lem:dimension_charts}
	Let $X$ be an algebraic space finitely presented over a field $k$.
	Suppose $a\colon U\to X$ is an \'{e}tale morphism from an affine scheme.
	\begin{enumerate}
		\item \label{lem:dimension_charts1}
		We have $\dim U\leq\dim X$, with equality if $a$ is surjective.
		\item \label{lem:dimension_charts2}
		If $u\in U$ maps to a closed point $x\in|X|$, then $u$ is closed in $U$ and $\dim(\mathcal{O}_{U,u})=\dim_x(X)$.
		\item \label{lem:dimension_charts3}
		Let $Z$ be an irreducible component of $U$.
		Suppose $W\subseteq U$ is an open subset containing the generic point of $Z$.
		Then there exists $u\in W\cap Z$ which lies in no other irreducible component of $U$, maps to a closed point of $|X|$, and satisfies $\dim(\mathcal{O}_{U,u})=\dim Z$.
	\end{enumerate}
\end{lemma}

\begin{proof}
	Since $X$ is finitely presented over $k$, it is quasi-separated, and hence, decent \cite[\href{https://stacks.math.columbia.edu/tag/03JX}{Tag 03JX}]{StacksProject}.
	As $a$ is \'{e}tale, $U$ is locally of finite type over $k$.
	Since $U$ is affine, it is of finite type over $k$, and hence, Noetherian.

	We prove \eqref{lem:dimension_charts1}.
	By \cite[\href{https://stacks.math.columbia.edu/tag/04N5}{Tag 04N5}]{StacksProject}, $\dim_u(U)=\dim_{a(u)}(X)$ for every $u\in U$.
	By \cite[\href{https://stacks.math.columbia.edu/tag/04MU}{Tag 04MU}]{StacksProject}, $\dim U$ is the supremum of $\dim_u(U)$ over $u\in U$.
	Hence, $\dim U$ is the supremum of $\dim_x(X)$ over the points $x$ in the image of $|a|$.
	This is at most $\dim X$, and equality occurs if $a$ is surjective.

	We prove \eqref{lem:dimension_charts2}.
	By \Cref{lem:closed_points_charts}, $u$ is closed in $U$.
	By \cite[\href{https://stacks.math.columbia.edu/tag/0A21}{Tag 0A21}]{StacksProject}, $\dim(\mathcal{O}_{U,u})=\dim_u(U)$.
	By \cite[\href{https://stacks.math.columbia.edu/tag/04N5}{Tag 04N5}]{StacksProject}, the latter is $\dim_x(X)$.

	We prove \eqref{lem:dimension_charts3}.
	Since $U$ is Noetherian, it has finitely many irreducible components \cite[\href{https://stacks.math.columbia.edu/tag/0052}{Tag 0052}]{StacksProject}.
	The generic point of $Z$ lies in no other irreducible component of $U$.
	Let $V\subseteq W$ be the complement in $W$ of the union of the irreducible components of $U$ other than $Z$.
	This union is closed, being a finite union of closed subsets, so $V$ is open in $U$.
	Moreover, $V$ contains the generic point of $Z$.
	Note that every point of $U$ lies in some irreducible component \cite[\href{https://stacks.math.columbia.edu/tag/004W}{Tag 004W}]{StacksProject}.
	Hence, $V\subseteq Z$, and each point of $V$ lies in no other irreducible component of $U$.
	Since $a$ is \'{e}tale, $|a|(V)$ is a nonempty open subset of $|X|$ \cite[\href{https://stacks.math.columbia.edu/tag/03BX}{Tag 03BX}]{StacksProject}.
	As $\operatorname{Spec}(k)$ is Jacobson, so is $X$ \cite[\href{https://stacks.math.columbia.edu/tag/0BA3}{Tag 0BA3}]{StacksProject}.
	Moreover, $X$ is decent, and so $|X|$ is Jacobson \cite[\href{https://stacks.math.columbia.edu/tag/0BA6}{Tag 0BA6}]{StacksProject}.
	Hence, $|X|$ is the closure of its set of closed points \cite[\href{https://stacks.math.columbia.edu/tag/005U}{Tag 005U}]{StacksProject}.
	In particular, the nonempty open subset $|a|(V)$ contains a closed point of $|X|$.
	Choose $u\in V$ mapping to such a closed point.
	By \eqref{lem:dimension_charts2}, $u$ is closed in $U$.
	As $u$ lies in no irreducible component of $U$ other than $Z$, \cite[\href{https://stacks.math.columbia.edu/tag/0A21}{Tag 0A21}]{StacksProject} gives
	\begin{displaymath}
		\dim(\mathcal{O}_{U,u})=\dim_u(U)=\dim Z.
	\end{displaymath}
	This completes the proof.
\end{proof}

\subsection{Support}
\label{sec:prelim_support}

For $M \in \operatorname{Qcoh}(Y)$, $\operatorname{Supp}(M) :=|p|(\operatorname{Supp}(p^\ast M)) \subseteq |Y|$ where $p \colon U \to Y$ is any \'{e}tale presentation from a scheme.
By \cite[\href{https://stacks.math.columbia.edu/tag/07TY}{Tag 07TY}]{StacksProject}, $\operatorname{Supp}$ is independent of the \'{e}tale presentation. 
It coincides with the notion of support for abelian sheaves on the small \'{e}tale site \cite[\href{https://stacks.math.columbia.edu/tag/04K7}{Tag 04K7}]{StacksProject}.
More generally, if $E \in D_{\operatorname{qc}}(Y)$, 
\begin{displaymath}
	\operatorname{Supp}(E) := \bigcup_{j \in \mathbb{Z}} \operatorname{Supp}(\mathcal{H}^j(E)) \subseteq |Y|.
\end{displaymath}
Given a closed subset $Z \subseteq |Y|$, $E\in D_{\operatorname{qc}}(Y)$ is \textbf{supported on $Z$} if $\operatorname{Supp}(E) \subseteq Z$. 
If $E\in \operatorname{Qcoh}(Y)$ is of finite type, then $\operatorname{Supp}(E)\subseteq |Y|$ is closed \cite[\href{https://stacks.math.columbia.edu/tag/07TZ}{Tag 07TZ}]{StacksProject}. 
Denote by $D_{\operatorname{qc},Z}(Y)$ the strictly full subcategory of $D_{\operatorname{qc}}(Y)$ consisting of objects supported on $Z$. 
If $Y$ is locally Noetherian, quasi-coherent sheaves of finite type coincide with coherent sheaves \cite[\href{https://stacks.math.columbia.edu/tag/07UB}{Tag 07UB}]{StacksProject}.
If $Y$ is Noetherian, then $D^b_{\operatorname{coh},Z}(Y) := D^b_{\operatorname{coh}}(Y) \cap D_{\operatorname{qc},Z}(Y)$.

\subsection{Duality}
\label{sec:prelim_duality}

We use the duality conventions of \cite[\S 2.4]{GuisadoVillaalgordo/Lank:2026} and \cite{Neeman:2023}.
Let $k$ be a field.
Consider a proper morphism $f\colon X\to S:=\operatorname{Spec}(k)$ of algebraic spaces.
Denote by $\epsilon_S\colon (S_{\textrm{et}},\mathcal{O}_{S_{\textrm{et}}}) \to (S_{\textrm{Zar}},\mathcal{O}_{S_{\textrm{Zar}}})$ the natural morphism of ringed topoi.
There exist adjoint equivalences:
\begin{equation}
	\label{eq:derived_equivalences_affine_topoi}
	\begin{aligned}
		\mathbf{R}(\epsilon_S)_{\ast}
		&\colon
		D_{\operatorname{qc}}(S_{\textrm{et}}) \leftrightarrows D_{\operatorname{qc}}(S_{\textrm{Zar}}) \colon \mathbf{L}\epsilon_S^{\ast},
		\\& \mathbf{R}\Gamma(S_{\textrm{Zar}},-)\colon
		D_{\operatorname{qc}}(S_{\textrm{Zar}}) \leftrightarrows D(k) \colon \widetilde{(-)}.
	\end{aligned}
\end{equation}
See \cite[\href{https://stacks.math.columbia.edu/tag/071Q}{Tags 071Q} and \href{https://stacks.math.columbia.edu/tag/06Z0}{06Z0}]{StacksProject}.
As $f$ is proper, $f^!$ is naturally isomorphic to the right adjoint $f^\times$ of $\mathbf{R}f_{\ast}$ on quasi-coherent derived categories.
See \cite[Lemma 2.22]{GuisadoVillaalgordo/Lank:2026}. 
Recall that $f^!$ is the upper shriek functor.
It is constructed from Nagata compactifications, and defined on $D_{\operatorname{qc}}$. 
See \cite[\S 2.4.2]{GuisadoVillaalgordo/Lank:2026}; the construction is carried out for algebraic spaces in \cite[Theorem 1.8]{Neeman:2023}.

\begin{lemma}
	[Serre duality]
	\label{lem:serre_duality}
	For every $K\in D_{\operatorname{qc}}(X)$, there exists a natural isomorphism
	\begin{equation}
		\label{eq:serre_duality}
		\operatorname{Hom}_X(K,f^!\mathcal{O}_{\operatorname{Spec}(k)})
		\cong\operatorname{Hom}_k(H^0(X_{\textrm{et}},K),k).
	\end{equation}
\end{lemma}

\begin{proof}
	By \cite[Lemma 2.22]{GuisadoVillaalgordo/Lank:2026},
	$f^!\mathcal{O}_{\operatorname{Spec}(k)}\cong f^\times\mathcal{O}_{\operatorname{Spec}(k)}$, where $f^\times$ is the right adjoint of $\mathbf{R}f_{\ast}$ on quasi-coherent derived categories.
	The adjunction between $\mathbf{R}f_{\ast}$ and $f^\times$ gives a natural isomorphism
	\begin{displaymath}
		\operatorname{Hom}_X(K,f^!\mathcal{O}_{\operatorname{Spec}(k)})
		\cong\operatorname{Hom}_{D_{\operatorname{qc}}(S_{\textrm{et}})}
		(\mathbf{R}f_{\ast}K,\mathcal{O}_{\operatorname{Spec}(k)}).
	\end{displaymath}
	The equivalences in \eqref{eq:derived_equivalences_affine_topoi} send
	$\mathcal{O}_{\operatorname{Spec}(k)}$ to $k$.
	Consequently, the group on the right is naturally isomorphic to
	\begin{displaymath}
		\operatorname{Hom}_{D(k)}(
		\mathbf{R}\Gamma(S_{\textrm{Zar}},
		\mathbf{R}(\epsilon_S)_{\ast}\mathbf{R}f_{\ast}K),k).
	\end{displaymath}
	By \cite[\href{https://stacks.math.columbia.edu/tag/0D6H}{Tag 0D6H}]{StacksProject},
	there is a natural isomorphism in $D(k)$,
	\begin{displaymath}
		\begin{aligned}
			&\mathbf{R}\Gamma(S_{\textrm{Zar}},
			\mathbf{R}(\epsilon_S)_{\ast}\mathbf{R}f_{\ast}K)\\
			&\qquad\cong
			\operatorname{Res}_{\Gamma(X_{\textrm{et}},\mathcal{O}_X)/k}
			\mathbf{R}\Gamma(X_{\textrm{et}},K).
		\end{aligned}
	\end{displaymath}
	Here the restriction-of-scalars functor is displayed to distinguish
	the derived category of $\Gamma(X_{\textrm{et}},\mathcal{O}_X)$-modules
	from $D(k)$.
	For every $A\in D(k)$, applying $H^0$ to morphisms yields a map
	\begin{displaymath}
		\operatorname{Hom}_{D(k)}(A,k)
		\to\operatorname{Hom}_k(H^0(A),k),
	\end{displaymath}
	which is natural in $A$ and involves no choices.
	It is bijective.
	Indeed, since $k$ is a field, $A$ is isomorphic in $D(k)$ to $\oplus_{n\in\mathbb{Z}}H^n(A)[-n]$.
	Under such an isomorphism, the source becomes the product of the groups $\operatorname{Hom}_{D(k)}(H^n(A)[-n],k)\cong\operatorname{Ext}^n_k(H^n(A),k)$.
	These vanish for $n\neq 0$, and the map becomes the projection onto the factor with $n=0$.
	Combining these isomorphisms proves \eqref{eq:serre_duality}.
\end{proof}

\begin{remark}
	\label{rmk:serre_duality_derived}
	It is possible to infer \Cref{lem:serre_duality} from \cite[\href{https://stacks.math.columbia.edu/tag/0E59}{Tag 0E59}]{StacksProject}.
\end{remark}

\begin{lemma}
	\label{lem:dualizing_support}
	$\operatorname{Supp}(f^!\mathcal{O}_{\operatorname{Spec}(k)})=|X|$.
\end{lemma}

\begin{proof}
	Fix $x\in|X|$.
	Choose an \'{e}tale germ $a\colon (U,u)\to (X,x)$ with $U$ an affine scheme.
	Denote by $f_U\colon U\to\operatorname{Spec}(k)$ the composite $f\circ a$.
	By \cite[Lemma 2.21]{GuisadoVillaalgordo/Lank:2026}, we have $(f\circ a)^!\cong a^!f^!$.
	As $U$ is affine and $X$ is quasi-separated, $a$ is quasi-affine \cite[Lemma 2.19]{GuisadoVillaalgordo/Lank:2026}, and in particular, separated.
	Hence, \cite[Lemma 3.2]{GuisadoVillaalgordo/Lank:2026b} applies to the separated \'{e}tale morphism $a$ of Noetherian algebraic spaces, and gives $a^!\cong\mathbf{L}a^{\ast}$.
	Consequently,
	\begin{displaymath}
		\mathbf{L}a^{\ast}f^!\mathcal{O}_{\operatorname{Spec}(k)}
		\cong a^!f^!\mathcal{O}_{\operatorname{Spec}(k)}
		\cong f_U^!\mathcal{O}_{\operatorname{Spec}(k)}.
	\end{displaymath}
	By \cite[\href{https://stacks.math.columbia.edu/tag/0AA3}{Tag 0AA3}]{StacksProject}, $f_U^!\mathcal{O}_{\operatorname{Spec}(k)}$ is a dualizing complex on $U$.
	Moreover, by \cite[\href{https://stacks.math.columbia.edu/tag/0A7G}{Tag 0A7G}]{StacksProject}, the stalk of $f_U^!\mathcal{O}_{\operatorname{Spec}(k)}$ at $u$ is a dualizing complex for $\mathcal{O}_{U,u}$.
	However, a dualizing complex $\omega$ for a nonzero ring $A$ is nonzero (e.g.\ \cite[\href{https://stacks.math.columbia.edu/tag/0A7B}{Tag 0A7B}]{StacksProject} implies $A\cong \mathbf{R}\operatorname{Hom}_A(\omega,\omega)$).
	Thus, $x\in\operatorname{Supp}(f^!\mathcal{O}_{\operatorname{Spec}(k)})$, which completes the proof.
\end{proof}

\subsection{Integral transforms}
\label{sec:prelim_integral_transforms}

Let $S$ be a Noetherian algebraic space.
Suppose $Y_1$ and $Y_2$ are proper and flat algebraic $S$-spaces.
Denote by $p_i\colon Y_1\times_S Y_2\to Y_i$ the canonical projections.
Then $Y_1$, $Y_2$, and $Y_1\times_SY_2$ are Noetherian, and the structure morphisms $Y_i\to S$ and the projections $p_i$ are proper, and so separated, quasi-compact, quasi-separated, and of finite presentation.
For $E\in D_{\operatorname{qc}}(Y_1\times_S Y_2)$, the \textbf{integral transform with kernel $E$} is the functor
\begin{displaymath}
	\Phi_E:=\mathbf{R}(p_2)_{\ast}(\mathbf{L}p_1^{\ast}(-)\otimes^{\mathbf{L}}E)\colon
	D_{\operatorname{qc}}(Y_1)\to D_{\operatorname{qc}}(Y_2).
\end{displaymath}
We say that $E\in D^b_{\operatorname{coh}}(Y_1\times_S Y_2)$ induces a \textbf{Fourier--Mukai equivalence} (over $S$) if $\Phi_E$ restricts to an equivalence $D^b_{\operatorname{coh}}(Y_1)\to D^b_{\operatorname{coh}}(Y_2)$.
In this case, $Y_1$ and $Y_2$ are called \textbf{Fourier--Mukai $S$-partners}, and we say it is \textbf{given} by $E$.
If clear from context, we omit the `$S$-' and say `Fourier--Mukai partners'.
We say that $A\in D_{\operatorname{qc}}(Y_1\times_S Y_2)$ is \textbf{relatively perfect over $Y_i$} if $A\otimes^{\mathbf{L}} \mathbf{L}p_i^\ast D^b_{\operatorname{qc}}(Y_i)\subseteq D^b_{\operatorname{qc}}(Y_1\times_S Y_2)$. 
By \cite[Corollary 1.4]{GuisadoVillaalgordo/Lank:2026}, an object $E\in D^-_{\operatorname{coh}}(Y_1\times_S Y_2)$ is relatively perfect over $Y_1$ (resp.\ $Y_2$) if, and only if, $\Phi_E (D^b_{\operatorname{coh}}(Y_1))\subseteq D^b_{\operatorname{coh}}(Y_2)$ (resp.\ $\Phi_E (\operatorname{Perf}(Y_1))\subseteq \operatorname{Perf}(Y_2)$).
In both cases, testing the defining condition on $\mathcal{O}_{Y_1}$ (resp.\ $\mathcal{O}_{Y_2}$) shows that such an $E$ has bounded cohomology.

\begin{lemma}
	\label{lem:transport_duality}
	Let $f_i \colon Y_i \to \operatorname{Spec}(k)$ be proper morphisms of algebraic spaces with $k$ a field.
	Consider a $k$-linear exact equivalence $F\colon D^b_{\operatorname{coh}}(Y_1)\to D^b_{\operatorname{coh}}(Y_2)$ which preserves perfect complexes.
	Then for every $P\in\operatorname{Perf}(Y_1)$ there exists a natural isomorphism,
	\begin{equation}\label{eq:transport_duality}
		F(P\otimes^{\mathbf{L}}f_1^!\mathcal{O}_{\operatorname{Spec}(k)})
		\cong F(P)\otimes^{\mathbf{L}}f_2^!\mathcal{O}_{\operatorname{Spec}(k)}.
	\end{equation}
\end{lemma}

\begin{proof}
	By hypothesis, $F(P)$ is perfect whenever $P$ is perfect. 
	By \cite[Lemma 3.3]{GuisadoVillaalgordo/Lank:2026b}, each $f_i^!\mathcal{O}_{\operatorname{Spec}(k)}$ is a dualizing complex on $Y_i$, and so it belongs to $D^b_{\operatorname{coh}}(Y_i)$ by definition.
	Hence, tensoring with a perfect complex gives a bounded pseudocoherent complex, so both sides of \eqref{eq:transport_duality} are defined and belong to $D^b_{\operatorname{coh}}(Y_2)$.
	We first recall Serre duality for $Y_1$.
	For $P\in\operatorname{Perf}(Y_1)$ and $Q\in D^b_{\operatorname{coh}}(Y_1)$, there exist isomorphisms which are natural in $P,Q$:
	\begin{equation}\label{eq:duality_hom}
		\begin{aligned}
			\operatorname{Hom}_{Y_1}
			& (Q,P\otimes^{\mathbf{L}}f_1^!\mathcal{O}_{\operatorname{Spec}(k)})
			\\&\cong\operatorname{Hom}_{Y_1}(Q,\operatorname{\mathbb{R}\mathcal{H}\! \mathit{om}}_{\mathcal{O}_{Y_1}}(P^\vee,f_1^!\mathcal{O}_{\operatorname{Spec}(k)})) && \textrm{\eqref{eq:perfect_dual}}
			\\&\cong\operatorname{Hom}_{Y_1}(Q\otimes^{\mathbf{L}}P^\vee,f_1^!\mathcal{O}_{\operatorname{Spec}(k)}) && (\textrm{\cite[\href{https://stacks.math.columbia.edu/tag/08J9}{Tag 08J9}]{StacksProject}})
			\\&\cong\operatorname{Hom}_k(H^0((Y_1)_{\textrm{et}},Q\otimes^{\mathbf{L}}P^\vee),k) && \textrm{(\Cref{lem:serre_duality})}
			\\&\cong\operatorname{Hom}_k(\operatorname{Hom}_{Y_1}(P,Q),k) && \textrm{\eqref{eq:perfect_dual}}.
		\end{aligned}
	\end{equation}
	Furthermore, for any $B\in D^b_{\operatorname{coh}}(Y_2)$, there exist isomorphisms natural in $B$ and $P$:
	\begin{displaymath}
		\begin{aligned}
			\operatorname{Hom}_{Y_2}(B,F(P\otimes^{\mathbf{L}}f_1^!\mathcal{O}_{\operatorname{Spec}(k)}))
			&\cong\operatorname{Hom}_{Y_1}(F^{-1}B,P\otimes^{\mathbf{L}}f_1^!\mathcal{O}_{\operatorname{Spec}(k)})
			\\&\cong\operatorname{Hom}_k(\operatorname{Hom}_{Y_1}(P,F^{-1}B),k) && \eqref{eq:duality_hom}
			\\&\cong\operatorname{Hom}_k(\operatorname{Hom}_{Y_2}(F(P),B),k)
			\\&\cong\operatorname{Hom}_{Y_2}(B,F(P)\otimes^{\mathbf{L}}f_2^!\mathcal{O}_{\operatorname{Spec}(k)}).
		\end{aligned}
	\end{displaymath}
	The first isomorphism uses a quasi-inverse $F^{-1}$ of $F$.
	The third is obtained by applying $\operatorname{Hom}_k(-,k)$ to the isomorphism $\operatorname{Hom}_{Y_1}(P,F^{-1}B)\cong\operatorname{Hom}_{Y_2}(F(P),B)$ induced by $F$. 
	The last is \eqref{eq:duality_hom} on $Y_2$.
	Both $F(P\otimes^{\mathbf{L}}f_1^!\mathcal{O}_{\operatorname{Spec}(k)})$ and $F(P)\otimes^{\mathbf{L}}f_2^!\mathcal{O}_{\operatorname{Spec}(k)}$ belong to $D^b_{\operatorname{coh}}(Y_2)$.
	Hence, by \cite[\href{https://stacks.math.columbia.edu/tag/001P}{Tag 001P}]{StacksProject}, there exists an isomorphism
	\begin{displaymath}
		F(P\otimes^{\mathbf{L}}f_1^!\mathcal{O}_{\operatorname{Spec}(k)}) \to F(P)\otimes^{\mathbf{L}}f_2^!\mathcal{O}_{\operatorname{Spec}(k)}
	\end{displaymath}
	which is natural in $P$.
	This completes the proof.
\end{proof}

\begin{remark}
	\Cref{lem:transport_duality} is the algebraic space analog of the Serre duality argument for schemes.
	See e.g.\ \cite[Theorem 25.127(3)]{Gortz/Wedhorn:2023}.
\end{remark}

\section{Local vanishing}
\label{sec:local_vanishing}

\begin{definition}
	\label{def:sup_for_abelian_bounded_above}
	Let $\mathcal{A}$ be an abelian category.
	For a nonzero $E \in D^-(\mathcal{A})$,
	\begin{displaymath}
		\operatorname{sup}(E):=\max\{i\in\mathbb{Z}: H^i(E)\neq 0\}.
	\end{displaymath} 
\end{definition}

\begin{lemma}
	\label{lem:top_tensor}
	Let $(R,\mathfrak{m})$ be a Noetherian local ring.
	Pick nonzero $A,B\in D^-_{\operatorname{coh}}(R)$.
	If $a=\operatorname{sup}(A)$ and $b=\operatorname{sup}(B)$, then $H^{a+b}(A\otimes^{\mathbf{L}}_R B)\cong H^a(A)\otimes_RH^b(B)\neq 0$ and $\operatorname{sup}(A\otimes^{\mathbf{L}}_RB)=a+b$.
\end{lemma}

\begin{proof}
	We prove the claim for the top degree of cohomological vanishing.
	Choose $U\in D^{\leq u}(R)$ and $V\in D^{\leq v}(R)$.
	Consider bounded above resolutions of $U$ and $V$ by flat modules with no nonzero terms respectively above degrees $u$ and $v$. 
	Denote these respectively by $U^\prime$ and $V^\prime$.
	See \cite[\href{https://stacks.math.columbia.edu/tag/064U}{Tag 064U}, \href{https://stacks.math.columbia.edu/tag/066E}{066E}, and \href{https://stacks.math.columbia.edu/tag/05T7}{05T7}]{StacksProject}.
	By \cite[\href{https://stacks.math.columbia.edu/tag/064K}{Tag 064K}]{StacksProject}, these representatives are $K$-flat.
	The associated total complex $\operatorname{Tot}(U^\prime \otimes V^\prime)$ has no terms in degrees larger than $u+v$.
	See \cite[\href{https://stacks.math.columbia.edu/tag/0GWN}{Tag 0GWN}]{StacksProject}.
	Hence, 
	\begin{equation}
		\label{eq:tensor_upper}
		D^{\leq u}(R)\otimes^{\mathbf{L}}_R D^{\leq v}(R)\subseteq D^{\leq u+v}(R).
	\end{equation}
	In particular, $U\otimes^{\mathbf{L}}_R V\in D^{\leq u+v}(R)$, and so $A\otimes^{\mathbf{L}}_R B\in D^{\leq a+b}(R)$.

	Next we prove the claimed canonical isomorphism.
	Since $H^i(A)=0$ for $i>a$, standard truncation triangles give
	\begin{displaymath}
		\tau^{\leq a-1}A\to A\to H^a(A)[-a]
		\to (\tau^{\leq a-1}A)[1].
	\end{displaymath}
	Tensoring this distinguished triangle with $B$ yields another distinguished triangle.
	By \eqref{eq:tensor_upper}, its left term is in $D^{\leq a+b-1}(R)$.
	By the long exact sequence in cohomology, 
	\begin{displaymath}
		\begin{aligned}
			0=H^{a+b}((\tau^{\leq a-1}A)\otimes^{\mathbf{L}} B)
			&\to H^{a+b}(A\otimes^{\mathbf{L}} B)\\
			&\to H^{a+b}(H^a(A)[-a]\otimes^{\mathbf{L}} B)\\
			&\to H^{a+b+1}((\tau^{\leq a-1}A)\otimes^{\mathbf{L}} B)=0.
		\end{aligned}
	\end{displaymath}
	Here both outer groups vanish by \eqref{eq:tensor_upper}.
	Indeed, $\tau^{\leq a-1}A\otimes^{\mathbf{L}}B\in D^{\leq a+b-1}(R)$.
	Thus, the middle morphism is an isomorphism.
	Moreover, by standard truncation triangles,
	\begin{displaymath}
		\tau^{\leq b-1}B\to B\to H^b(B)[-b]
		\to (\tau^{\leq b-1}B)[1].
	\end{displaymath}
	After tensoring it with $H^a(A)[-a]$, \eqref{eq:tensor_upper} implies the left term belongs to $D^{\leq a+b-1}(R)$.
	The long exact sequence in cohomology yields
	\begin{displaymath}
		\begin{aligned}
			H^{a+b}(H^a(A)[-a]\otimes^{\mathbf{L}} B)
			&\cong H^{a+b}(H^a(A)[-a]\otimes^{\mathbf{L}} H^b(B)[-b])
			\\&=H^0(H^a(A)\otimes^{\mathbf{L}} H^b(B))
			\\&=H^a(A)\otimes H^b(B).
		\end{aligned}
	\end{displaymath}
	The final equality uses 
	\begin{displaymath}
		H^0(H^a(A)\otimes^{\mathbf{L}} H^b(B))=\operatorname{Tor}_0^R(H^a(A),H^b(B)).
	\end{displaymath}
	The desired canonical isomorphism follows.

	We verify nonvanishing.
	Since $H^a(A)$ and $H^b(B)$ are finite nonzero modules over the local ring $(R,\mathfrak{m})$, \cite[\href{https://stacks.math.columbia.edu/tag/07RC}{Tag 07RC}]{StacksProject} implies $H^a(A)/\mathfrak{m} H^a(A) \neq 0$ and $H^b(B)/\mathfrak{m} H^b(B)\neq 0$.
	Set $\kappa=R/\mathfrak{m}$.
	Associativity of tensor product gives
	\begin{displaymath}
		(H^a(A) \otimes H^b(B))\otimes \kappa
		\cong(H^a(A) \otimes \kappa)\otimes_\kappa( H^b(B) \otimes\kappa).
	\end{displaymath}
	The right-hand side is the tensor product of two nonzero $\kappa$-vector spaces, and so it is nonzero.
	Hence, $H^a(A)\otimes H^b(B)\neq 0$.
	Consequently, $H^{a+b}(A\otimes^{\mathbf{L}} B)\neq 0$, proving $\operatorname{sup}(A\otimes^{\mathbf{L}} B)=a+b$.
\end{proof}

\begin{lemma}
	\label{lem:local-depth-dimension}
	Let $X$ be a locally Noetherian algebraic space.
	Choose $x\in|X|$.
	If $a\colon U\to X$ is an \'{e}tale morphism from a scheme and $u\in U$ maps to $x$, then
	$\operatorname{depth}(\mathcal{O}_{U,u})$ and $\dim(\mathcal{O}_{U,u})$
	are independent of the pair $(U,u)$.
\end{lemma}

\begin{proof}
	By \cite[\href{https://stacks.math.columbia.edu/tag/0BAM}{Tag 0BAM}]{StacksProject}, independence of the Krull dimension of the local ring follows.
	We prove the claim for depth.
	Fix an integer $d\geq0$.
	Consider the property of a germ $(U,u)\to (X,x)$ that $\operatorname{depth}(\mathcal{O}_{U,u})=d$.
	By shrinking, we may assume the germs are affine, and hence, Noetherian.
	If $(V,v)\to(U,u)$ is an \'{e}tale morphism of Noetherian scheme germs,
	the induced local homomorphism $\mathcal{O}_{U,u}\to\mathcal{O}_{V,v}$ is \'{e}tale \cite[\href{https://stacks.math.columbia.edu/tag/039N}{Tag 039N}]{StacksProject}.
	By \cite[\href{https://stacks.math.columbia.edu/tag/039T}{Tag 039T}]{StacksProject},
	the two local rings have the same depth.
	Hence, the property is \'{e}tale local in the sense of
	\cite[\href{https://stacks.math.columbia.edu/tag/04N1}{Tag 04N1}]{StacksProject}.
	Applying \cite[\href{https://stacks.math.columbia.edu/tag/04N2}{Tag 04N2}]{StacksProject}, the value of $d$ is independent of the chosen \'{e}tale germ about $x$.
	This completes the proof.
\end{proof}

\begin{definition}
	\label{def:local-depth-dimension}
	With notation as in \Cref{lem:local-depth-dimension}, the
	\textbf{depth} and the \textbf{Krull dimension of the local ring of $X$ at $x$}
	are respectively the integers $\operatorname{depth}(\mathcal{O}_{U,u})$ and
	$\dim(\mathcal{O}_{U,u})$.
	These are respectively denoted by $\operatorname{depth}_x(X)$ and
	$\dim_x^{\mathrm{loc}}(X)$.
\end{definition}

\begin{lemma}
	\label{lem:local-sup-space}
	Let $g\colon X^\prime\to X$ be a flat morphism of quasi-separated locally Noetherian algebraic spaces.
	Choose $x^\prime\in|X^\prime|$ with image $x\in|X|$.
	Consider \'{e}tale germs $a\colon(U,u)\to(X,x)$ and $a^\prime\colon(U^\prime,u^\prime)\to(X^\prime,x^\prime)$.
	Pick $K\in D^-_{\operatorname{coh}}(X)$ and $i\in\mathbb{Z}$.
	Then $\mathcal{H}^i(\mathbf{L}a^{\ast}K)_u=0$ if, and only if, $\mathcal{H}^i(\mathbf{L}(a^\prime)^{\ast}\mathbf{L}g^{\ast}K)_{u^\prime}=0$. 
\end{lemma}

\begin{proof}
	By shrinking, we may assume $U$ and $U^\prime$ are affine schemes. 
	Form the fiber product of $a$ and $g\circ a^\prime$:
	\begin{displaymath}
		\begin{tikzcd}[column sep=large,row sep=large]
			U^\prime\times_X U \arrow[r] \arrow[d] & U \arrow[d,"a"]\\
			U^\prime \arrow[r,"g\circ a^\prime"'] & X.
		\end{tikzcd}
	\end{displaymath}
	By \cite[Lemma 2.19]{GuisadoVillaalgordo/Lank:2026}, $a$ and $a^\prime$ are quasi-affine, and hence representable.
	Hence, $U^\prime\times_X U$ is a scheme.
	By \cite[\href{https://stacks.math.columbia.edu/tag/03H4}{Tag 03H4}]{StacksProject}, there exists $z\in U^\prime\times_X U$ mapping to both $u^\prime$ and $u$.
	Base change implies the canonical projections $U^\prime\times_X U\to U$ and $U^\prime\times_X U\to U^\prime$ are flat.
	By \cite[\href{https://stacks.math.columbia.edu/tag/00HR}{Tag 00HR}]{StacksProject}, the induced local homomorphisms $\mathcal{O}_{U,u}\to\mathcal{O}_{U^\prime\times_X U,z}$ and $\mathcal{O}_{U^\prime,u^\prime}\to\mathcal{O}_{U^\prime\times_X U,z}$ are faithfully flat.
	Denote by $p$ and $q$ respectively their induced morphisms on affine spectra.
	There exists a commutative diagram 
	\begin{displaymath}
		\begin{tikzcd}
			{\operatorname{Spec}(\mathcal{O}_{U^\prime\times_X U,z})} && {\operatorname{Spec}(\mathcal{O}_{U,u})} \\
			& {U^\prime\times_X U} & U \\
			{\operatorname{Spec}(\mathcal{O}_{U^\prime,u^\prime})} & {U^\prime} & X
			\arrow["p", from=1-1, to=1-3]
			\arrow[from=1-1, to=2-2]
			\arrow["q"', from=1-1, to=3-1]
			\arrow["{f_1}", from=1-3, to=2-3]
			\arrow[from=2-2, to=2-3]
			\arrow[from=2-2, to=3-2]
			\arrow["a", from=2-3, to=3-3]
			\arrow["{f_2}"', from=3-1, to=3-2]
			\arrow["{g\circ a^\prime}"', from=3-2, to=3-3]
		\end{tikzcd}
	\end{displaymath}
	where $f_i$ are the canonical morphisms coming from localization.
	Since the diagram commutes, $a\circ f_1\circ p=g\circ a^\prime\circ f_2\circ q$.
	Recall that pullback along flat morphisms commutes with cohomology sheaves.
	Consequently, for every $i\in\mathbb{Z}$,
	\begin{displaymath}
		\begin{aligned}
			p^\ast f_1^\ast\mathcal{H}^i(\mathbf{L}a^{\ast}K)
			&\cong\mathcal{H}^i(\mathbf{L}(a\circ f_1 \circ p)^{\ast}K)
			\\&\cong\mathcal{H}^i(\mathbf{L}(g\circ a^\prime \circ f_2 \circ q)^{\ast}K)
			\\&\cong q^\ast f_2^\ast\mathcal{H}^i(\mathbf{L}(a^\prime)^{\ast}\mathbf{L}g^{\ast}K).
		\end{aligned}
	\end{displaymath}
	Here $f_1^\ast\mathcal{H}^i(\mathbf{L}a^{\ast}K)$ is the quasi-coherent sheaf on $\operatorname{Spec}(\mathcal{O}_{U,u})$ associated to the stalk $\mathcal{H}^i(\mathbf{L}a^{\ast}K)_u$, and similarly, $f_2^\ast\mathcal{H}^i(\mathbf{L}(a^\prime)^{\ast}\mathbf{L}g^{\ast}K)$ corresponds to the stalk at $u^\prime$.
	Since $p$ and $q$ are faithfully flat, pullback along them reflects the vanishing of quasi-coherent sheaves.
	See e.g.\ \cite[Lemma 12.81]{Gortz/Wedhorn:2020}.
	Therefore, $\mathcal{H}^i(\mathbf{L}a^{\ast}K)_u=0$ if, and only if, $\mathcal{H}^i(\mathbf{L}(a^\prime)^{\ast}\mathbf{L}g^{\ast}K)_{u^\prime}=0$.
\end{proof}

\begin{remark-definition}
	\label{def:pointwise_sup}
	Let $X$ be a locally Noetherian quasi-separated algebraic space.
	Choose $K\in D^-_{\operatorname{coh}}(X)$ and $x\in|X|$.
	By \Cref{lem:local-sup-space} applied to the identity of $X$, the set
	\begin{displaymath}
		\{i\in\mathbb{Z}:\mathcal{H}^i(\mathbf{L}a^{\ast}K)_u\neq 0\},
	\end{displaymath}
	does not depend on the choice of \'{e}tale germ $a\colon(U,u)\to(X,x)$.
	The \textbf{support} $\operatorname{Supp}(K)$ consists of the points
	for which this set is nonempty (which agrees with that of \Cref{sec:prelim_support}).
	If $x\in\operatorname{Supp}(K)$, then the set of integers of nonvanishing of cohomology is bounded above.
	Define
	\begin{displaymath}
		\operatorname{sup}_x(K)
		:=\max\{i\in\mathbb{Z}:\mathcal{H}^i(\mathbf{L}a^{\ast}K)_u\neq0\}.
	\end{displaymath}
	We call $\operatorname{sup}_x(K)$ the \textbf{pointwise supremum} of $K$ at $x$.
\end{remark-definition}

\begin{corollary}
	\label{cor:flat_pullback_sup}
	Let $X$ be a Noetherian algebraic space and $K\in D^-_{\operatorname{coh}}(X)$.
	Suppose $g\colon X^\prime\to X$ is a flat morphism from a Noetherian algebraic space. 
	Choose $x^\prime\in|X^\prime|$ which maps to $x\in|X|$.
	Then $x^\prime\in\operatorname{Supp}(\mathbf{L}g^{\ast}K)$ if, and only if, $x\in\operatorname{Supp}(K)$.
	If these conditions hold, then
	\begin{displaymath}
		\operatorname{sup}_{x^\prime}(\mathbf{L}g^{\ast}K)
		=\operatorname{sup}_x(K).
	\end{displaymath}
	Moreover, if $K\in D^b_{\operatorname{coh}}(X)$ is nonzero, then
	\begin{displaymath}
		\operatorname{sup}(K)
		=\max\{\operatorname{sup}_x(K):x\in\operatorname{Supp}(K)\}.
	\end{displaymath}
	(Note that $D^-_{\operatorname{coh}}(X) \subseteq D^-(\mathcal{O}_X)$, and so we are applying \Cref{def:sup_for_abelian_bounded_above}.)
\end{corollary}

\begin{proof}
	Since $g$ is flat, $\mathbf{L}g^{\ast}K$ has bounded above coherent cohomology.
	Applying \Cref{lem:local-sup-space} gives the equivalence of support
	conditions and equality of the corresponding sets of nonzero cohomological degrees.
	Taking their maximum gives the desired equality.

	We check the last claim.
	A quasi-coherent sheaf on $X$ vanishes if, and only if, its support is empty.
	Thus, $\operatorname{sup}(K)$ is the largest integer $i$ for which
	$\mathcal{H}^i(K)$ has nonempty support.
	By the definition of the pointwise supremum, this is precisely
	the maximum in the statement.
\end{proof}

\begin{lemma}
	\label{lem:pointwise_tensor}
	Let $X$ be a Noetherian algebraic space.
	Choose $K,L\in D^-_{\operatorname{coh}}(X)$.
	Then $K\otimes^{\mathbf{L}}_{\mathcal{O}_X}L \in D^-_{\operatorname{coh}}(X)$, and
	\begin{displaymath}
		\operatorname{Supp}(K\otimes^{\mathbf{L}}_{\mathcal{O}_X}L)
		=\operatorname{Supp}(K)\cap\operatorname{Supp}(L).
	\end{displaymath}
	Moreover, for every point $x$ in this intersection,
	\begin{displaymath}
		\operatorname{sup}_x(K\otimes^{\mathbf{L}}_{\mathcal{O}_X}L)
		=\operatorname{sup}_x(K)+\operatorname{sup}_x(L).
	\end{displaymath}
\end{lemma}

\begin{proof}
	By \cite[\href{https://stacks.math.columbia.edu/tag/08IK}{Tags 08IK} and \href{https://stacks.math.columbia.edu/tag/09J9}{09J9}]{StacksProject}, the objects of $D^-_{\operatorname{coh}}(X)$ are exactly the pseudocoherent complexes, and $D^-_{\operatorname{coh}}(X)$ is closed under tensor product.
	Hence, $K\otimes^{\mathbf{L}}_{\mathcal{O}_X}L\in D^-_{\operatorname{coh}}(X)$.
	Fix $x\in|X|$.
	Choose an \'{e}tale germ $a\colon(U,u)\to(X,x)$.
	There exists an isomorphism,
	\begin{displaymath}
		(\mathbf{L}a^{\ast}(K\otimes^{\mathbf{L}}_{\mathcal{O}_X}L))_u
		\cong (\mathbf{L}a^{\ast}K)_u \otimes^{\mathbf{L}}_{\mathcal{O}_{U,u}} (\mathbf{L}a^{\ast}L)_u.
	\end{displaymath}
	Both factors on the right belong to
	$D^-_{\operatorname{coh}}(\mathcal{O}_{U,u})$.
	If at least one factor is zero, then their derived tensor product is zero.
	If their derived tensor product is zero, then \Cref{lem:top_tensor} implies at least one of the factors is zero.
	Hence, the desired intersection holds. 
	Moreover, by \Cref{lem:top_tensor} for $x$ in the intersection, the tensor product is nonzero and has top cohomology in degree
	$\operatorname{sup}_x(K)+\operatorname{sup}_x(L)$.
	This proves both assertions.
\end{proof}

\section{Points in supports}
\label{sec:point_supported}

\begin{lemma}
	\label{lem:point_perfect}
	Let $X$ be a Noetherian algebraic space.
	\begin{enumerate}
		\item \label{lem:point_perfect1}
		For every closed point $x\in|X|$, there exists
		$P_x\in\operatorname{Perf}(X)$ such that
		$\operatorname{Supp}(P_x)=\{x\}$.
		\item \label{lem:point_perfect2}
		Let $x\in|X|$ be a closed point.
		Choose $P_x\in\operatorname{Perf}(X)$ satisfying
		$\operatorname{Supp}(P_x)=\{x\}$.
		Then, for every $Q\in D^b_{\operatorname{coh}}(X)$,
		\begin{displaymath}
			x\in\operatorname{Supp}(Q)
			\quad\Longleftrightarrow\quad
			\operatorname{Hom}_X(P_x,Q[i])\neq 0
			\textrm{ for some }i\in\mathbb{Z}.
		\end{displaymath}
		\item \label{lem:point_perfect3}
		Choose $P_x$ as in \eqref{lem:point_perfect2} for every closed
		point $x\in|X|$.
		If $Q\in D^b_{\operatorname{coh}}(X)$ is nonzero, then
		$\operatorname{Hom}_X(P_x,Q[i])\neq 0$ for some closed point
		$x\in|X|$ and some $i\in\mathbb{Z}$.
	\end{enumerate}
\end{lemma}

\begin{proof}
	We prove \eqref{lem:point_perfect1}.
	Since $|X|$ is a Noetherian topological space \cite[\href{https://stacks.math.columbia.edu/tag/04ZF}{Tag 04ZF}]{StacksProject}, the open subset $|X|\setminus\{x\}$ is quasi-compact.
		Hence, as $X$ is quasi-compact and quasi-separated, \cite[\href{https://stacks.math.columbia.edu/tag/0AEC}{Tag 0AEC}]{StacksProject} provides a perfect object $P_x$ which generates $D_{\mathrm{qc},\{x\}}(X)$ and is supported on $\{x\}$.
	It remains to see that this support is nonempty.
	Since $X$ is decent, there is a closed residual immersion $i_x\colon\operatorname{Spec}(\kappa(x))\to X$ \cite[\href{https://stacks.math.columbia.edu/tag/03JX}{Tags 03JX}, \href{https://stacks.math.columbia.edu/tag/0H1R}{0H1R}, and \href{https://stacks.math.columbia.edu/tag/0H1U}{0H1U}]{StacksProject}.
	The sheaf $(i_x)_{\ast}\mathcal{O}_{\operatorname{Spec}(\kappa(x))}$
	is nonzero, coherent, and supported on $\{x\}$.
	Hence, $D_{\mathrm{qc},\{x\}}(X)$ is nonzero.
	A compact generator of a nonzero category is nonzero, and a nonzero object has nonempty support.
	Thus, $\operatorname{Supp}(P_x)=\{x\}$.

	Before proving \eqref{lem:point_perfect2}, we establish an observation. 
	Let $M$ be a coherent sheaf on $X$ with $\operatorname{Supp}(M)\subseteq\{x\}$.
	We claim that
	\begin{equation}
		\label{eq:closed_point_sheaf_acyclicity}
		H^j(X_{\textrm{et}},M)=0\quad(j>0),
		\textrm{ and }
		M\neq 0\ \Longrightarrow\ \Gamma(X_{\textrm{et}},M)\neq 0.
	\end{equation}
	Denote by $\mathcal{I}\subseteq\mathcal{O}_X$ the ideal sheaf defining the
	closed residual immersion $i_x$ (again, $x\in |X|$ is a closed point).
	Since $M$ is coherent and supported set-theoretically at $x$,
	there exists an integer $n\geq1$ such that $\mathcal{I}^nM=0$.
	See \cite[\href{https://stacks.math.columbia.edu/tag/07UK}{Tag 07UK}]{StacksProject}.
	Consider the finite filtration
	\begin{displaymath}
		M\supseteq\mathcal{I}M\supseteq\cdots
		\supseteq\mathcal{I}^{n-1}M
		\supseteq\mathcal{I}^nM=0.
	\end{displaymath}
	Note that each quotient $\mathcal{I}^jM/\mathcal{I}^{j+1}M$
	is coherent and annihilated by $\mathcal{I}$.
	By \cite[\href{https://stacks.math.columbia.edu/tag/08AM}{Tag 08AM}]{StacksProject}, there exists a finite-dimensional $\kappa(x)$-vector space $V_j$ such that $\mathcal{I}^jM/\mathcal{I}^{j+1}M \cong(i_x)_{\ast}\widetilde{V_j}$.
	Since $i_x$ is finite, \cite[\href{https://stacks.math.columbia.edu/tag/0A4K}{Tag 0A4K}]{StacksProject} implies its pushforward on the small \'{e}tale site is exact.
	Using derived global sections and the \'{e}tale--Zariski cohomology comparison, we obtain
	\begin{displaymath}
		\begin{aligned}
			\mathbf{R}\Gamma(X_{\textrm{et}},(i_x)_{\ast}\widetilde{V_j})
			&\cong\mathbf{R}\Gamma(
				(\operatorname{Spec}(\kappa(x)))_{\textrm{et}},
				\widetilde{V_j})\\
			&\cong V_j[0].
		\end{aligned}
	\end{displaymath}
	See \cite[\href{https://stacks.math.columbia.edu/tag/0D6H}{Tags 0D6H},
	\href{https://stacks.math.columbia.edu/tag/03DW}{03DW}, and
	\href{https://stacks.math.columbia.edu/tag/01XB}{01XB}]{StacksProject}.
	Since $\mathcal{I}^0:=\mathcal{O}_X$, it follows that $\mathcal{I}^0M=M$. 
	Then the short exact sequences
	\begin{displaymath}
		0\to\mathcal{I}^{j+1}M\to\mathcal{I}^jM
		\to(i_x)_{\ast}\widetilde{V_j}\to0
	\end{displaymath}
	and descending induction on $j$ show that $H^i(X_{\textrm{et}},\mathcal{I}^jM)=0$ for every $i>0$ and every $j\in\{0,\ldots,n\}$.
	Moreover, global sections yield a surjection
	\begin{displaymath}
		\Gamma(X_{\textrm{et}},M)\to V_0.
	\end{displaymath}
	If $M\neq 0$, then $V_0\neq 0$.
	Indeed, otherwise $M=\mathcal{I}M$, and a further iteration above gives $M=\mathcal{I}^n M=0$.
	This proves \eqref{eq:closed_point_sheaf_acyclicity}.

	Now we prove \eqref{lem:point_perfect2}.
	Let $Q\in D^b_{\operatorname{coh}}(X)$.
	Set $C=P_x^\vee\otimes^{\mathbf{L}}_{\mathcal{O}_X}Q$.
	Since $P_x$ is perfect, $P_x^\vee$ is perfect \cite[\href{https://stacks.math.columbia.edu/tag/08JJ}{Tag 08JJ}]{StacksProject}, and so both have closed support \cite[\href{https://stacks.math.columbia.edu/tag/07TZ}{Tag 07TZ}]{StacksProject}.
	As derived internal Homs commute with restriction to open subspaces \cite[\href{https://stacks.math.columbia.edu/tag/08JB}{Tag 08JB}]{StacksProject} and $(P_x^\vee)^\vee\cong P_x$, the complexes $P_x$ and $P_x^\vee$ vanish on the same open subspaces of $X$.
	Hence, $\operatorname{Supp}(P_x^\vee)=\operatorname{Supp}(P_x)=\{x\}$.
	Hence, $C\in D^b_{\operatorname{coh}}(X)$, and \Cref{lem:pointwise_tensor} gives $\operatorname{Supp}(C)=\{x\}\cap\operatorname{Supp}(Q)$.
	In particular, $C\neq 0$ if, and only if, $x\in\operatorname{Supp}(Q)$.
	Note that every cohomology sheaf of $C$ is coherent and supported on $\{x\}$.
	We claim that for every $j\in\mathbb{Z}$ there exists a natural isomorphism
	\begin{equation}
		\label{eq:closed_point_acyclicity}
		H^j(X_{\textrm{et}},C)
		\cong\Gamma(X_{\textrm{et}},\mathcal{H}^j(C)).
	\end{equation}
	By induction on the number of nonzero cohomology sheaves of $C$, this follows from \eqref{eq:closed_point_sheaf_acyclicity}.
	Indeed, if $b=\operatorname{sup}(C)$, apply derived global sections to the standard truncation triangle
	\begin{displaymath}
		\tau^{\leq b-1}C\to C\to\mathcal{H}^b(C)[-b] \to ( \tau^{\leq b-1}C)[1].
	\end{displaymath}
		Since $\mathcal{H}^b(C)$ is a coherent sheaf supported on $\{x\}$, \eqref{eq:closed_point_sheaf_acyclicity} applies to it, and gives
	\begin{displaymath}
		H^j(X_{\textrm{et}},\mathcal{H}^b(C))=0
	\end{displaymath}
	for $j>0$.
	Since a sheaf sits in cohomological degree $0$, these groups also vanish for $j<0$.
	Hence,
	\begin{displaymath}
		\mathbf{R}\Gamma(X_{\textrm{et}},\mathcal{H}^b(C))
		\cong\Gamma(X_{\textrm{et}},\mathcal{H}^b(C))[0].
	\end{displaymath}
	As $\mathbf{R}\Gamma(X_{\textrm{et}},-)$ commutes with shifts,
	\begin{displaymath}
		\mathbf{R}\Gamma(X_{\textrm{et}},\mathcal{H}^b(C)[-b])
		\cong\Gamma(X_{\textrm{et}},\mathcal{H}^b(C))[-b].
	\end{displaymath}
	In particular, the derived global sections of $\mathcal{H}^b(C)[-b]$ are concentrated in degree $b$, with cohomology $\Gamma(X_{\textrm{et}},\mathcal{H}^b(C))$ in that degree.
	The induction hypothesis, applied to $\tau^{\leq b-1}C$, and the long exact cohomology sequence of the resulting distinguished triangle give \eqref{eq:closed_point_acyclicity} in degrees below $b$ and in degree $b$; in degrees above $b$ both sides vanish, because $\tau^{\leq b-1}C$ has no cohomology in degree $\geq b$ and $\mathcal{H}^j(C)=0$ for $j>b$.
	The case $C=0$ is immediate.
	For every $j\in\mathbb{Z}$, \eqref{eq:perfect_dual} and \eqref{eq:closed_point_acyclicity} yield
	\begin{displaymath}
		\begin{aligned}
			\operatorname{Hom}_X(P_x,Q[j])
			&\cong H^0(X_{\textrm{et}},C[j])\\
			&=H^j(X_{\textrm{et}},C)\\
			&\cong\Gamma(X_{\textrm{et}},\mathcal{H}^j(C)).
		\end{aligned}
	\end{displaymath}
	By \eqref{eq:closed_point_sheaf_acyclicity}, the last group is nonzero if, and only if, $\mathcal{H}^j(C)\neq 0$.
	Thus, $\operatorname{Hom}_X(P_x,Q[j])\neq 0$ for some $j$ if, and only if, $C\neq 0$.
	The latter is equivalent to $x\in\operatorname{Supp}(Q)$.

	Finally, we prove \eqref{lem:point_perfect3}.
	If $Q\neq 0$, one of its coherent cohomology sheaves has nonempty closed support $Z\subseteq|X|$.
	Since $X$ is decent, $|X|$ is Kolmogorov \cite[\href{https://stacks.math.columbia.edu/tag/03K6}{Tag 03K6}]{StacksProject}.
	As $Z$ is nonempty, quasi-compact, and Kolmogorov,
	\cite[\href{https://stacks.math.columbia.edu/tag/005E}{Tag 005E}]{StacksProject} gives a point $x$ closed in $Z$, and hence in $|X|$.
	Applying \eqref{lem:point_perfect2} at this point finishes the proof.
\end{proof}

\begin{lemma}
	\label{lem:rouquier_thickening_spaces}
	Let $X$ be a Noetherian algebraic space.
	Suppose $i\colon Z\to X$ is a closed immersion defined by an ideal sheaf $\mathcal{I}\subseteq\mathcal{O}_X$.
	For every integer $n\geq1$, let $i_n\colon Z_n=V(\mathcal{I}^n)\to X$ denote the corresponding closed immersion.
	Then for every $C\in D^b_{\operatorname{coh},|i|(|Z|)}(X)$ there exist an integer $n\geq1$ and an object $C_n\in D^b_{\operatorname{coh}}(Z_n)$ such that $C\cong\mathbf{R}(i_n)_{\ast}C_n$ in $D^b_{\operatorname{coh}}(X)$.
\end{lemma}

\begin{proof}
	We adapt the proof of \cite[Lemma 7.40]{Rouquier:2008} and \cite[Remark 2.23]{Lank:2024} to Noetherian algebraic spaces.
	It holds nearly verbatim but we include the details.
	By \cite[Lemma 4.7]{Lank:2026}, we may represent $C$ by a bounded complex $C^\bullet$ of coherent sheaves on $X$.
	By \cite[\href{https://stacks.math.columbia.edu/tag/07UK}{Tag 07UK}]{StacksProject}, for every coherent sheaf $M$ on $X$ with $\operatorname{Supp}(M)\subseteq|Z|$, there exists an integer $d\geq1$ such that $\mathcal{I}^d M=0$.
	We claim that $C^\bullet$ is quasi-isomorphic to a bounded complex of coherent sheaves whose terms are supported on $|Z|$.
	We prove this by induction on the number of terms of $C^\bullet$ not supported on $|Z|$.
	If every term is supported on $|Z|$, there is nothing to prove.
	Otherwise, choose the smallest integer $i$ such that $\operatorname{Supp}(C^i)\nsubseteq|Z|$.
	Denote by $d^i\colon C^i\to C^{i+1}$ the differential.
	There exists an exact sequence
	\begin{displaymath}
		0\to\operatorname{im}(d^{i-1})
		\to\ker(d^i)
		\to\mathcal{H}^i(C^\bullet)\to0.
	\end{displaymath}
	Since $C^{i-1}$ and $\mathcal{H}^i(C^\bullet)$ are supported on $|Z|$, the same holds for $\ker(d^i)$.
	By the preceding observation, there exists an integer $d\geq1$ such that $\mathcal{I}^d\ker(d^i)=0$. 
	By \cite[\href{https://stacks.math.columbia.edu/tag/07UL}{Tag 07UL}]{StacksProject}, there exists an integer $c\geq0$ such that for every $m\geq c$,
	\begin{displaymath}
		\mathcal{I}^mC^i\cap\ker(d^i)
		=
		\mathcal{I}^{m-c}
		(\mathcal{I}^cC^i\cap\ker(d^i)).
	\end{displaymath}
	Choose an integer $n\geq\max\{1,c+d\}$.
	Since $\mathcal{I}^d\ker(d^i)=0$, it follows that $\mathcal{I}^n C^i \cap\ker(d^i)=0$.
	Hence, the restriction of $d^i$ induces an isomorphism $\mathcal{I}^n C^i \to d^i(\mathcal{I}^n C^i)$.
	Let $R^\bullet\subseteq C^\bullet$ be the subcomplex consisting of these two sheaves in degrees $i$ and $i+1$, and zero in all other degrees.
	Since the differential is an isomorphism, $R^\bullet$ is contractible.
	Consequently, the natural morphism $C^\bullet\to C^\bullet/R^\bullet$ is a quasi-isomorphism.
	The term of the quotient in degree $i$ is $C^i/\mathcal{I}^nC^i$, which is annihilated by $\mathcal{I}^n$, and hence supported on $|Z|$.
	Note that the preceding terms are unchanged.
	Moreover, the remaining terms are unchanged except in degree $i+1$, where the term is a quotient of $C^{i+1}$.
	Thus, $C^\bullet/R^\bullet$ has strictly fewer terms not supported on $|Z|$.
	Applying the induction hypothesis proves the claim.
	Hence, we can assume that every term $C^j$ is supported on $|Z|$.
	Since $C^\bullet$ is bounded, there exists an integer $N\geq1$ such that $\mathcal{I}^N C^j =0$ for $j\in\mathbb{Z}$.
	By \cite[\href{https://stacks.math.columbia.edu/tag/08AM}{Tag 08AM}]{StacksProject}, each $C^j$ is the pushforward of a coherent sheaf
	$C_N^j$ on $Z_N$.
	Moreover, the pushforward functor identifies coherent sheaves on $Z_N$ with coherent sheaves on $X$ annihilated by $\mathcal{I}^N$.
	Since this equivalence is fully faithful, the differentials of $C^\bullet$ lift uniquely to morphisms $d_N^j\colon C_N^j\to C_N^{j+1}$.
	These morphisms define a bounded complex $C_N^\bullet$ of coherent sheaves on $Z_N$, and there is an isomorphism $(i_N)_{\ast}C_N^\bullet\cong C^\bullet$.
	Since $i_N$ is finite, its pushforward on the small \'{e}tale site is exact
	\cite[\href{https://stacks.math.columbia.edu/tag/0A4K}{Tag 0A4K}]{StacksProject}.
	Hence, $\mathbf{R}(i_N)_{\ast}C_N^\bullet \cong(i_N)_{\ast}C_N^\bullet \cong C^\bullet$.
	Taking $C_N\in D^b_{\operatorname{coh}}(Z_N)$ to be the object represented by $C_N^\bullet$ proves the claim.
\end{proof}

\begin{lemma}
	\label{lem:finite_support_exact}
	Let $q\colon V\to T$ be a proper morphism of Noetherian algebraic spaces.
	Suppose $p\colon W\to V$ is a closed immersion such that $q \circ p$ is finite.
	Then for every $C\in D^b_{\operatorname{coh},|p|(|W|)}(V)$ there exist isomorphisms $\mathcal{H}^j(\mathbf{R}q_{\ast}C) \cong q_{\ast}\mathcal{H}^j(C)$ for every $j\in\mathbb{Z}$.
	Moreover, if $M$ is a coherent sheaf on $V$ supported on $|W|$, then $q_{\ast}M=0$ if, and only if, $M=0$.
	Consequently, for every nonzero $C\in D^b_{\operatorname{coh},|p|(|W|)}(V)$, $\operatorname{sup}(\mathbf{R}q_{\ast}C) = \operatorname{sup}(C)$.
\end{lemma}

\begin{proof}
	Let $\mathcal{I}\subseteq\mathcal{O}_V$ be the ideal sheaf associated to $p$.
	By \Cref{lem:rouquier_thickening_spaces}, there exist an integer $n\geq1$ and $C_n\in D^b_{\operatorname{coh}}(W_n)$ such that $C\cong\mathbf{R}(i_n)_{\ast}C_n$ and $i_n\colon W_n=V(\mathcal{I}^n)\to V$.
	Since $q$ is proper, $q\circ i_n\colon W_n\to T$ is proper.
	Let $\operatorname{Spec}(K)\to T$ be a morphism where $K$ is a field.
	The closed immersion $W\to W_n$ is a thickening, and hence, so is its base change along $W_n\times_T\operatorname{Spec}(K)\to W_n$ \cite[\href{https://stacks.math.columbia.edu/tag/09ZX}{Tag 09ZX}]{StacksProject}.
	Thus, the fibers of $q\circ p$ and $q\circ i_n$ over $\operatorname{Spec}(K)$ have the same underlying topological spaces.
	As $q\circ p$ is finite, these are discrete \cite[\href{https://stacks.math.columbia.edu/tag/0A4X}{Tag 0A4X}]{StacksProject}.
	Hence, by loc.\ cit., the proper morphism $q\circ i_n$ is finite.
	By \cite[\href{https://stacks.math.columbia.edu/tag/0A4K}{Tag 0A4K}]{StacksProject}, the pushforward along a finite morphism is exact
	on the small \'{e}tale site.
	Hence, $\mathbf{R}(q\circ i_n)_\ast$ commutes with cohomology sheaves.
	Moreover, a finite morphism is quasi-affine, and so \cite[Lemma 3.2]{GuisadoVillaalgordo/Lank:2026} ensures that $\mathbf{R}(q\circ i_n)_\ast$ is conservative (i.e.\ reflects the zero objects).
	Since $\mathbf{R}q_{\ast}C\cong\mathbf{R}(q\circ i_n)_{\ast}C_n$ and $\mathcal{H}^j(C)\cong(i_n)_{\ast}\mathcal{H}^j(C_n)$, it follows that
	\begin{displaymath}
		\mathcal{H}^j(\mathbf{R}q_{\ast}C)
		\cong(q\circ i_n)_{\ast}\mathcal{H}^j(C_n)
		\cong q_{\ast}\mathcal{H}^j(C)
	\end{displaymath}
	for every $j\in\mathbb{Z}$.
	Similarly, if $M$ is a coherent sheaf on $V$ supported on $|W|$ and $M\cong(i_n)_{\ast}M_n$, then $q_{\ast}M\cong(q\circ i_n)_{\ast}M_n$, which is nonzero whenever $M_n$, and hence $M$, is nonzero.
	Consequently, $\operatorname{sup}(\mathbf{R}q_{\ast}C)=\operatorname{sup}(C)$ for every nonzero $C$.
\end{proof}

\section{Depth at closed points}
\label{sec:depth}

\begin{lemma}
	\label{lem:point_tensor_comparison}
	Consider proper morphisms of algebraic spaces $f_i \colon Y_i \to \operatorname{Spec}(k)$ where $k$ is a field.
	Let $x\in|Y_1|$ be a closed point, $y\in|Y_2|$, and $E\in D^b_{\operatorname{coh}}(Y_1\times_k Y_2)$ be relatively perfect over $Y_1$.
	Choose $P_x\in\operatorname{Perf}(Y_1)$ such that $\operatorname{Supp}(P_x)=\{x\}$.
	Assume $y\in\operatorname{Supp}(\Phi_E(P_x))$.
	Then $y\in\operatorname{Supp}(\Phi_E(P_x\otimes^{\mathbf{L}}f_1^!\mathcal{O}_{\operatorname{Spec}(k)}))$, and
	\begin{equation}\label{eq:point_tensor_comparison}
		\operatorname{sup}_y(\Phi_E(P_x\otimes^{\mathbf{L}}
		f_1^!\mathcal{O}_{\operatorname{Spec}(k)}))
		=\operatorname{sup}_y(\Phi_E(P_x)) +\operatorname{sup}_x(f_1^!\mathcal{O}_{\operatorname{Spec}(k)}).
	\end{equation}
\end{lemma}

\begin{proof}
	Choose an \'{e}tale germ $a\colon (U,u)\to (Y_2,y)$.
	Denote by $h\colon\operatorname{Spec}(\mathcal{O}_{U,u})\to Y_2$ the composite of $a$ with the canonical morphism $\operatorname{Spec}(\mathcal{O}_{U,u})\to U$, and by $\mathfrak{m}_u$ the closed point of $\operatorname{Spec}(\mathcal{O}_{U,u})$.
	Let $p_i \colon Y_1\times_k Y_2 \to Y_i$ be the canonical projections.
	Since localization is flat and $a$ is \'{e}tale, $h$ is flat.
	Consider the fibered square
	\begin{displaymath}
		\begin{tikzcd}[column sep=large,row sep=large]
			Y_1\times_k\operatorname{Spec}(\mathcal{O}_{U,u})
				\arrow[r,"g"]\arrow[d,"q"']
			&Y_1\times_kY_2\arrow[d,"p_2"]\\
			\operatorname{Spec}(\mathcal{O}_{U,u})
				\arrow[r,"h"'] &Y_2.
		\end{tikzcd}
	\end{displaymath}
	Write $p=p_1\circ g$ for the first projection.
	Since both $h$ and the structure morphism $\operatorname{Spec}(\mathcal{O}_{U,u})\to\operatorname{Spec}(k)$ are flat, base change implies that $g$ and $p$ are flat.
	Moreover, by base change, $f_1$ being proper implies $q$ is as well.
	Note that $Y_1\times_k\operatorname{Spec}(\mathcal{O}_{U,u})$ is Noetherian 
	because it is finitely presented over the Noetherian scheme $\operatorname{Spec}(\mathcal{O}_{U,u})$.

	Since $h$ is flat, the square is Tor-independent.
	By \cite[\href{https://stacks.math.columbia.edu/tag/08IR}{Tags 08IR} and \href{https://stacks.math.columbia.edu/tag/07A4}{07A4}]{StacksProject}, there exist natural isomorphisms for every $A\in D_{\operatorname{qc}}(Y_1)$,
	\begin{equation}
		\label{eq:point_comparison_flat_base_change}
		\mathbf{L}h^{\ast}\Phi_E(A)
		\cong\mathbf{R}q_{\ast}\mathbf{L}g^{\ast}
			(\mathbf{L}p_1^{\ast}A\otimes^{\mathbf{L}}E)
		\cong\mathbf{R}q_{\ast}(
			\mathbf{L}p^{\ast}A\otimes^{\mathbf{L}}
			\mathbf{L}g^{\ast}E).
	\end{equation}
	Here we used that $p_1\circ g=p$.
		Since every module $N$ over the local ring $\mathcal{O}_{U,u}$ satisfies $N_{\mathfrak{m}_u}=N$, a quasi-coherent sheaf on $\operatorname{Spec}(\mathcal{O}_{U,u})$ vanishes if, and only if, its stalk at $\mathfrak{m}_u$ vanishes.
	Therefore, for every $A\in D^-_{\operatorname{coh}}(Y_2)$, \Cref{cor:flat_pullback_sup} gives
	\begin{equation}
		\label{eq:point_comparison_local}
		y\in\operatorname{Supp}(A)
		\quad\Longleftrightarrow\quad
		\mathbf{L}h^{\ast}A\neq0,
	\end{equation}
	and, when these equivalent conditions hold,
	\begin{displaymath}
		\operatorname{sup}_y(A)=\operatorname{sup}(\mathbf{L}h^{\ast}A).
	\end{displaymath}

	Denote by $i_x\colon\operatorname{Spec}(\kappa(x))\to Y_1$ the closed residual immersion associated to $x$.
	As $f_1$ is proper with affine target, $Y_1$ is separated, and hence, decent \cite[\href{https://stacks.math.columbia.edu/tag/03JX}{Tag 03JX}]{StacksProject}.
	Moreover, $x$ is closed, so \cite[\href{https://stacks.math.columbia.edu/tag/0H1R}{Tags 0H1R} and \href{https://stacks.math.columbia.edu/tag/0H1U}{0H1U}]{StacksProject} implies $i_x$ is a closed immersion.
	In particular, $f_1\circ i_x$ is finitely presented, so \cite[\href{https://stacks.math.columbia.edu/tag/00FV}{Tag 00FV}]{StacksProject} says $\kappa(x)/k$ is finite.
	The base change of $i_x$ along $p$ is the closed immersion
	\begin{displaymath}
		\operatorname{Spec}(\kappa(x))\times_k
		\operatorname{Spec}(\mathcal{O}_{U,u})
		\to Y_1\times_k\operatorname{Spec}(\mathcal{O}_{U,u}).
	\end{displaymath}
	By \cite[\href{https://stacks.math.columbia.edu/tag/03H4}{Tag 03H4}]{StacksProject}, the underlying topological space of $\operatorname{Spec}(\kappa(x))\times_k \operatorname{Spec}(\mathcal{O}_{U,u})$ is $p^{-1}(\{x\})$.
	Moreover, $\operatorname{Spec}(\kappa(x))\times_k \operatorname{Spec}(\mathcal{O}_{U,u})$ is finite over $\operatorname{Spec}(\mathcal{O}_{U,u})$ because it is the base change of the finite morphism $\operatorname{Spec}(\kappa(x))\to\operatorname{Spec}(k)$.

	Put $C=\mathbf{L}p^{\ast}P_x\otimes^{\mathbf{L}}\mathbf{L}g^{\ast}E$.
	Since $p$ is flat and $\operatorname{Supp}(P_x)=\{x\}$,
	\Cref{cor:flat_pullback_sup} gives
	\begin{displaymath}
		\operatorname{Supp}(\mathbf{L}p^{\ast}P_x)
		=|\operatorname{Spec}(\kappa(x))\times_k
		\operatorname{Spec}(\mathcal{O}_{U,u})|.
	\end{displaymath}
	By \Cref{lem:pointwise_tensor}, the support of $C$
	is contained in this closed subset.
	Since $g$ is flat, $\mathbf{L}g^{\ast}E$ has bounded pseudocoherent cohomology.
	Furthermore, $\mathbf{L}p^{\ast}P_x$ is perfect, so it follows that $C\in D^b_{\operatorname{coh}}(Y_1\times_k\operatorname{Spec}(\mathcal{O}_{U,u}))$.
	Moreover, by \cite[Lemma 3.3]{GuisadoVillaalgordo/Lank:2026b}, $f_1^!\mathcal{O}_{\operatorname{Spec}(k)}\in D^b_{\operatorname{coh}}(Y_1)$.
	Hence, $P_x\otimes^{\mathbf{L}}f_1^!\mathcal{O}_{\operatorname{Spec}(k)}$ is bounded pseudocoherent.
	As $E$ is relatively perfect over $Y_1$, it follows that
	\begin{displaymath}
		\mathbf{L}p_1^{\ast}(P_x\otimes^{\mathbf{L}}f_1^!\mathcal{O}_{\operatorname{Spec}(k)}) \otimes^{\mathbf{L}}E \in D^b_{\operatorname{coh}}(Y_1\times_kY_2).
	\end{displaymath}
	Furthermore, by \cite[Corollary 1.4]{GuisadoVillaalgordo/Lank:2026}, relative perfectness of $E$ over $Y_1$ implies that $\Phi_E$ maps $D^b_{\operatorname{coh}}(Y_1)$ into $D^b_{\operatorname{coh}}(Y_2)$.
	In particular, both $\Phi_E(P_x)$ and $\Phi_E(P_x\otimes^{\mathbf{L}}f_1^!\mathcal{O}_{\operatorname{Spec}(k)})$ belong to $D^b_{\operatorname{coh}}(Y_2)$.
	Flat pullback along $g$ with \cite[\href{https://stacks.math.columbia.edu/tag/07A4}{Tag 07A4}]{StacksProject} gives
	\begin{displaymath}
		\begin{aligned}
			&C\otimes^{\mathbf{L}}
				\mathbf{L}p^{\ast}(f_1^!\mathcal{O}_{\operatorname{Spec}(k)})\\
			&\quad\cong
				\mathbf{L}g^{\ast}
				(\mathbf{L}p_1^{\ast}
				(P_x\otimes^{\mathbf{L}}f_1^!\mathcal{O}_{\operatorname{Spec}(k)})
				\otimes^{\mathbf{L}}E)
				\in D^b_{\operatorname{coh}}(Y_1\times_k\operatorname{Spec}(\mathcal{O}_{U,u})).
		\end{aligned}
	\end{displaymath}
	By \Cref{lem:pointwise_tensor}, its support is contained in
	$\operatorname{Supp}(C)$, and hence in
	$|\operatorname{Spec}(\kappa(x))\times_k
	\operatorname{Spec}(\mathcal{O}_{U,u})|$.
	Thus, \Cref{lem:finite_support_exact} applies to $q$, the displayed
	closed immersion, and both complexes $C$ and $C\otimes^{\mathbf{L}}\mathbf{L}p^{\ast}(f_1^!\mathcal{O}_{\operatorname{Spec}(k)})$.
	Indeed, each of them is nonzero if, and only if, its derived pushforward along $q$ is nonzero, and in such a case, $\mathbf{R}q_{\ast}$ preserves its supremum.

	We first treat $\Phi_E(P_x)\in D^b_{\operatorname{coh}}(Y_2)$.
	By \eqref{eq:point_comparison_flat_base_change}, $\mathbf{L}h^{\ast}\Phi_E(P_x)\cong\mathbf{R}q_{\ast}C$.
	Since $y\in\operatorname{Supp}(\Phi_E(P_x))$, \eqref{eq:point_comparison_local} shows that $\mathbf{R}q_{\ast}C\neq0$.
	Hence, $C\neq0$ and
	\begin{equation}
		\label{eq:point_comparison_first_degree}
		\operatorname{sup}_y(\Phi_E(P_x))
		=\operatorname{sup}(\mathbf{R}q_{\ast}C)
		=\operatorname{sup}(C).
	\end{equation}

	Next let $n:=\operatorname{sup}_x(f_1^!\mathcal{O}_{\operatorname{Spec}(k)})$.
	Note that $n$ is defined by \Cref{lem:dualizing_support}.
	Let $w\in\operatorname{Supp}(C)$.
	Since $\operatorname{Supp}(C)\subseteq p^{-1}(\{x\})$, $w$ maps to $x$.
	Hence, \Cref{cor:flat_pullback_sup} gives $w\in\operatorname{Supp}(\mathbf{L}p^{\ast}(f_1^!\mathcal{O}_{\operatorname{Spec}(k)}))$ and $\operatorname{sup}_w(\mathbf{L}p^{\ast}(f_1^!\mathcal{O}_{\operatorname{Spec}(k)}))=n$.
	By \Cref{lem:pointwise_tensor}, it follows that
	\begin{displaymath}
		\operatorname{Supp}(C\otimes^{\mathbf{L}}
			\mathbf{L}p^{\ast}(f_1^!\mathcal{O}_{\operatorname{Spec}(k)}))
		=\operatorname{Supp}(C),
	\end{displaymath}
	as well as for every $w\in\operatorname{Supp}(C)$,
	\begin{displaymath}
		\operatorname{sup}_w(C\otimes^{\mathbf{L}}
			\mathbf{L}p^{\ast}(f_1^!\mathcal{O}_{\operatorname{Spec}(k)}))
		=\operatorname{sup}_w(C)+n.
	\end{displaymath}
	Recall that $C\neq0$.
	Taking the maximum over $w\in\operatorname{Supp}(C)$, and using \Cref{cor:flat_pullback_sup}, it follows that $C\otimes^{\mathbf{L}}\mathbf{L}p^{\ast}(f_1^!\mathcal{O}_{\operatorname{Spec}(k)})$ is nonzero and
	\begin{equation}
		\label{eq:point_comparison_degree_shift}
		\operatorname{sup}(C\otimes^{\mathbf{L}}
			\mathbf{L}p^{\ast}(f_1^!\mathcal{O}_{\operatorname{Spec}(k)}))
		=\operatorname{sup}(C)+n.
	\end{equation}

	Finally, by \eqref{eq:point_comparison_flat_base_change},
	\begin{displaymath}
		\mathbf{L}h^{\ast}\Phi_E(P_x\otimes^{\mathbf{L}}
			f_1^!\mathcal{O}_{\operatorname{Spec}(k)})
		\cong\mathbf{R}q_{\ast}(C\otimes^{\mathbf{L}}
			\mathbf{L}p^{\ast}(f_1^!\mathcal{O}_{\operatorname{Spec}(k)})).
	\end{displaymath}
	By \Cref{lem:finite_support_exact} and \eqref{eq:point_comparison_degree_shift}, the complex above is nonzero with supremum $\operatorname{sup}(C)+n$.
	Since $\Phi_E(P_x\otimes^{\mathbf{L}}f_1^!\mathcal{O}_{\operatorname{Spec}(k)})\in D^b_{\operatorname{coh}}(Y_2)$, \eqref{eq:point_comparison_local} shows that $y\in\operatorname{Supp}(\Phi_E(P_x\otimes^{\mathbf{L}}f_1^!\mathcal{O}_{\operatorname{Spec}(k)}))$ and
	\begin{displaymath}
		\operatorname{sup}_y(\Phi_E(P_x\otimes^{\mathbf{L}}
			f_1^!\mathcal{O}_{\operatorname{Spec}(k)}))
		=\operatorname{sup}(C)+n.
	\end{displaymath}
	Combined with \eqref{eq:point_comparison_first_degree}, this proves \eqref{eq:point_tensor_comparison}.
\end{proof}

\begin{theorem}
	\label{thm:depth_correspondence}
	Let $f_i \colon Y_i \to \operatorname{Spec}(k)$ be proper morphisms of algebraic spaces with $k$ a field. 
	Suppose $Y_1$ and $Y_2$ are Fourier--Mukai $k$-partners given by $E\in D^b_{\operatorname{coh}}(Y_1\times_k  Y_2)$.
	For any closed point $x\in|Y_1|$, pick $P_x\in\operatorname{Perf}(Y_1)$ with $\operatorname{Supp}(P_x)=\{x\}$.
	Then for every closed point $y\in\operatorname{Supp}(\Phi_E(P_x))$,
	\begin{equation}
		\label{eq:depth_correspondence}
		\operatorname{depth}_x(Y_1)=\operatorname{depth}_y(Y_2).
	\end{equation}
	Furthermore, for every closed point $x\in|Y_1|$, the support of $\Phi_E(P_x)$ contains a closed point, every closed point of $Y_2$ lies in $\operatorname{Supp}(\Phi_E(P_x))$ for some closed point $x\in|Y_1|$, and
	\begin{equation}\label{eq:depth_spectra}
		\{\operatorname{depth}_x(Y_1):x\in |Y_1|\textrm{ closed}\}
		=\{\operatorname{depth}_y(Y_2):y\in |Y_2|\textrm{ closed}\}.
	\end{equation}
\end{theorem}

\begin{proof}
	By \cite[Remark 2.7]{GuisadoVillaalgordo/Lank:2026}, $E$ is relatively perfect over both projections.
	Hence, \cite[Corollary 1.4]{GuisadoVillaalgordo/Lank:2026} shows that $\Phi_E$ preserves perfect complexes.
	Since $\Phi_E$ is an integral transform over $k$, it is $k$-linear.
	Applying \Cref{lem:transport_duality} gives
	\begin{equation}
		\label{eq:duality_spaces}
		\Phi_E(P_x\otimes^{\mathbf{L}}f_1^!\mathcal{O}_{\operatorname{Spec}(k)})
		\cong\Phi_E(P_x)\otimes^{\mathbf{L}}f_2^!\mathcal{O}_{\operatorname{Spec}(k)}.
	\end{equation}

	We identify the top degrees of the dualizing complexes at closed points.
	Choose an \'{e}tale germ $a\colon (U,u)\to (Y_1,x)$ from an affine scheme.
	Since $Y_1$ is separated, it is decent \cite[\href{https://stacks.math.columbia.edu/tag/03JX}{Tag 03JX}]{StacksProject}. 
	Hence, by \Cref{lem:closed_points_charts}, $u$ is closed.
	Set $f_U=f_1\circ a$. 
	By \cite[Lemma 3.2]{GuisadoVillaalgordo/Lank:2026b} and \cite[Lemma 2.21]{GuisadoVillaalgordo/Lank:2026},
	\begin{displaymath}
		\mathbf{L}a^{\ast}f_1^!\mathcal{O}_{\operatorname{Spec}(k)}
		\cong f_U^!\mathcal{O}_{\operatorname{Spec}(k)}.
	\end{displaymath}
	By \cite[\href{https://stacks.math.columbia.edu/tag/0AWL}{Tag 0AWL}]{StacksProject}, the right-hand side is a dualizing complex normalized at the closed point $u$.
	The same argument applies to $Y_2$.
		Since a normalized dualizing complex has vanishing cohomology in positive degrees, \cite[\href{https://stacks.math.columbia.edu/tag/0A7U}{Tags 0A7U} and \href{https://stacks.math.columbia.edu/tag/0AWR}{0AWR}]{StacksProject} implies
	\begin{equation}
		\label{eq:depth_normalized_spaces}
		\begin{aligned}
			\operatorname{sup}_x(f_1^!\mathcal{O}_{\operatorname{Spec}(k)})
			&=-\operatorname{depth}_x(Y_1),\\
			\operatorname{sup}_y(f_2^!\mathcal{O}_{\operatorname{Spec}(k)})
			&=-\operatorname{depth}_y(Y_2).
		\end{aligned}
	\end{equation}

	Fix a closed point $x\in|Y_1|$ and a closed point $y\in\operatorname{Supp}(\Phi_E(P_x))$.
	Put $s=\operatorname{sup}_y(\Phi_E(P_x))$.
	By \Cref{lem:point_tensor_comparison} and \eqref{eq:depth_normalized_spaces}, the left-hand side of \eqref{eq:duality_spaces} has top degree $s-\operatorname{depth}_x(Y_1)$ at $y$.
	By \Cref{lem:pointwise_tensor} and \Cref{lem:dualizing_support}, the right-hand side has top degree $s-\operatorname{depth}_y(Y_2)$.
	These two complexes are isomorphic, and so the two degrees agree.
	This proves \eqref{eq:depth_correspondence}.

	We next show that every closed point of $Y_2$ occurs in the support of some $\Phi_E(P_x)$.
	Choose a closed point $y\in|Y_2|$.
	Denote by $i_y\colon\operatorname{Spec}(\kappa(y))\to Y_2$ its associated closed residual immersion.
	Note that $(i_y)_{\ast}\mathcal{O}_{\operatorname{Spec}(\kappa(y))}$
	is a nonzero coherent sheaf.
	Set $Q=\Phi_E^{-1}((i_y)_{\ast}\mathcal{O}_{\operatorname{Spec}(\kappa(y))})$.
	By \Cref{lem:point_perfect}, there exist a closed point $x\in|Y_1|$ and an integer $i$ such that $\operatorname{Hom}_{Y_1}(P_x,Q[i])\neq 0$.
	The equivalence gives
	\begin{displaymath}
		\operatorname{Hom}_{Y_2}(
		\Phi_E(P_x),(i_y)_{\ast}
		\mathcal{O}_{\operatorname{Spec}(\kappa(y))}[i])\neq0.
	\end{displaymath}
	Since $\Phi_E(P_x)$ is perfect, \eqref{eq:perfect_dual} identifies this Hom-set with
	\begin{displaymath}
		\begin{aligned}
			H^i((Y_2)_{\textrm{et}},
			\Phi_E(P_x)^\vee\otimes^{\mathbf{L}}
			(i_y)_{\ast}\mathcal{O}_{\operatorname{Spec}(\kappa(y))}).
		\end{aligned}
	\end{displaymath}
	If $y\notin\operatorname{Supp}(\Phi_E(P_x))$, then
	$\Phi_E(P_x)^\vee$ and
	$(i_y)_{\ast}\mathcal{O}_{\operatorname{Spec}(\kappa(y))}$
	have disjoint supports (see the proof of \Cref{lem:point_perfect}\eqref{lem:point_perfect2}).
	Hence, by \Cref{lem:pointwise_tensor}, their derived tensor product vanishes, which is a contradiction.
	Consequently, it follows that $y\in\operatorname{Supp}(\Phi_E(P_x))$.

		Conversely, for every closed point $x\in|Y_1|$, $P_x$ is nonzero.
	Hence, $\Phi_E(P_x)\in D^b_{\operatorname{coh}}(Y_2)$ is nonzero, so some cohomology sheaf $\mathcal{H}^j(\Phi_E(P_x))$ is nonzero, and its support is a nonempty closed subset of $|Y_2|$ contained in $\operatorname{Supp}(\Phi_E(P_x))$ \cite[\href{https://stacks.math.columbia.edu/tag/07TZ}{Tag 07TZ}]{StacksProject}.
	Since $Y_2$ is decent and of finite type over $k$, $|Y_2|$ is Jacobson \cite[\href{https://stacks.math.columbia.edu/tag/0BA3}{Tags 0BA3} and \href{https://stacks.math.columbia.edu/tag/0BA6}{0BA6}]{StacksProject}, so this closed subset contains a point which is closed in $|Y_2|$.
	Hence, $\operatorname{Supp}(\Phi_E(P_x))$ contains a closed point $y$.
		Combined with the preceding paragraph, this proves \eqref{eq:depth_spectra}.
\end{proof}

\begin{example}
	\label{exm:poincare}
	Let $k$ be an algebraically closed field.
	Suppose $A$ is an elliptic curve over $k$.
	Identify $A$ with its dual abelian variety via the principal polarization given by the origin.
	Denote by $\mathcal{P}$ the normalized Poincar\'{e} bundle on $A\times_kA$.
	By \cite{Mukai:1981}, $\Phi_{\mathcal{P}}\colon D^b_{\operatorname{coh}}(A)\to D^b_{\operatorname{coh}}(A)$ is an equivalence; see also \cite[Chapter 9]{Huybrechts:2006}.
	For a closed point $a\in A$, the skyscraper sheaf $\mathcal{O}_a$ is perfect.
	Indeed, $A$ is regular, and $\operatorname{Supp}(\mathcal{O}_a)=\{a\}$.
	Moreover, its transform $\Phi_{\mathcal{P}}(\mathcal{O}_a)$ is the restriction of $\mathcal{P}$ to $\{a\}\times_kA\cong A$, which is an invertible sheaf.
	In particular, $\operatorname{Supp}(\Phi_{\mathcal{P}}(\mathcal{O}_a))=|A|$.
	Here every closed point of $A$ has depth $1$ (cf.\ \eqref{eq:depth_correspondence}).
	Thus, the support appearing in \Cref{thm:depth_correspondence} need not be finite.
\end{example}

\section{Invariants}
\label{sec:dimension}

\begin{definition}
	\label{def:cohen_macaulay}
	Let $X$ be a locally Noetherian algebraic space.
	We say that $X$ is \textbf{Cohen--Macaulay at} $x\in|X|$ if $\operatorname{depth}_x(X)=\dim_x^{\mathrm{loc}}(X)$.
	If this holds at every point of $|X|$, we say $X$ is \textbf{Cohen--Macaulay}.
\end{definition}

\begin{remark}
	By \Cref{lem:local-depth-dimension}, $X$ is Cohen--Macaulay at $x$ if, and only if, $\mathcal{O}_{U,u}$ is Cohen--Macaulay for some (equivalently, every) \'{e}tale germ $(U,u)\to(X,x)$ with $U$ a scheme.
	Thus, \Cref{def:cohen_macaulay} coincides with the notion from \cite[\href{https://stacks.math.columbia.edu/tag/04RC}{Tags 04RC} and \href{https://stacks.math.columbia.edu/tag/0BBL}{0BBL}]{StacksProject}.
\end{remark}

\begin{lemma}
	\label{lem:cm_closed_points}
	Let $X$ be an algebraic space which is finitely presented over a field $k$.
	Then the set of points at which $X$ is Cohen--Macaulay is open in $|X|$.
	In particular, $X$ is Cohen--Macaulay if, and only if,
	it is Cohen--Macaulay at every closed point of $|X|$.
\end{lemma}

\begin{proof}
	Choose an \'{e}tale presentation $U\to X$ with $U$ an affine
	scheme.
	By \Cref{lem:local-depth-dimension}, the inverse image in $U$ of
	the Cohen--Macaulay locus of $X$ is the Cohen--Macaulay locus of $U$.
	By \cite[\href{https://stacks.math.columbia.edu/tag/00RF}{Tag 00RF}]{StacksProject}, the Cohen--Macaulay locus of $U$ is open.
	Since $U\to X$ is open and surjective, the Cohen--Macaulay locus of
	$X$ is open as well.

	We prove the last assertion.
	By \cite[\href{https://stacks.math.columbia.edu/tag/03JX}{Tag 03JX}]{StacksProject}, $X$ is decent.
	Moreover, by
	\cite[\href{https://stacks.math.columbia.edu/tag/0BA3}{Tags 0BA3} and \href{https://stacks.math.columbia.edu/tag/0BA6}{0BA6}]{StacksProject}, $|X|$ is Jacobson.
	If the complement of the Cohen--Macaulay locus were nonempty, it would contain a closed point.
	Hence, Cohen--Macaulayness at every closed point implies
	Cohen--Macaulayness everywhere.
	The converse is immediate.
\end{proof}

\begin{lemma}
	\label{lem:max_closed_depth}
	Let $X$ be a nonempty algebraic space which is finitely presented over a field $k$.
	Then
	\begin{displaymath}
		\dim X
		=\max\{\operatorname{depth}_x(X):x\in|X|\textrm{ closed}\}.
	\end{displaymath}
\end{lemma}

\begin{proof}
	Set $d=\dim X$.
	Choose an \'{e}tale presentation $a\colon U\to X$ with $U$ an affine scheme.
	By \Cref{lem:dimension_charts}\eqref{lem:dimension_charts1}, $\dim U=d$.
	As $U$ is Noetherian, it has finitely many irreducible components \cite[\href{https://stacks.math.columbia.edu/tag/0052}{Tag 0052}]{StacksProject}.
	Note that every chain of irreducible closed subsets of $U$ is contained in one of them. 
	Indeed, its largest member is contained in an irreducible component \cite[\href{https://stacks.math.columbia.edu/tag/004W}{Tag 004W}]{StacksProject}.
	Hence, there exists an irreducible component $Z$ of $U$ of Krull dimension $d$.

	Let $W\subseteq U$ be the set of points at which the local ring of $U$ is Cohen--Macaulay.
	By \cite[\href{https://stacks.math.columbia.edu/tag/00RF}{Tag 00RF}]{StacksProject}, $W$ is open.
	Since the generic point of $Z$ corresponds to a minimal prime ideal, its local ring has Krull dimension zero, and hence, must be Cohen--Macaulay.
	Thus, $W$ contains the generic point of $Z$.
	By \Cref{lem:dimension_charts}\eqref{lem:dimension_charts3}, there exists $u\in W\cap Z$ which maps to a closed point $x\in|X|$ and satisfies $\dim(\mathcal{O}_{U,u})=\dim Z=d$.
	Since $\mathcal{O}_{U,u}$ is Cohen--Macaulay, its depth is $d$.
	By \Cref{def:local-depth-dimension}, $\operatorname{depth}_x(X)=d$.

	Conversely, let $y\in|X|$ be a closed point.
	Choose an \'{e}tale germ $b\colon(U^\prime,u^\prime)\to(X,y)$ with $U^\prime$ an affine scheme.
	By \cite[\href{https://stacks.math.columbia.edu/tag/00LK}{Tag 00LK} and \href{https://stacks.math.columbia.edu/tag/04MU}{04MU}]{StacksProject} and \Cref{lem:dimension_charts}\eqref{lem:dimension_charts1},
	\begin{displaymath}
		\operatorname{depth}_y(X)
		=\operatorname{depth}(\mathcal{O}_{U^\prime,u^\prime})
		\leq\dim(\mathcal{O}_{U^\prime,u^\prime})
		\leq\dim U^\prime
		\leq d.
	\end{displaymath}
		Consequently, $d$ is the maximum of the depths at closed points.
\end{proof}

\begin{corollary}
	\label{cor:dimension_invariance}
	Let $f_i \colon Y_i \to \operatorname{Spec}(k)$ be proper morphisms of algebraic spaces with $k$ a field. 
	Assume that $Y_1$ and $Y_2$ are Fourier--Mukai partners given by some $E\in D^b_{\operatorname{coh}}(Y_1\times_k Y_2)$. 
	Then $\dim Y_1=\dim Y_2$.
\end{corollary}

\begin{proof}
	If $Y_1$ is empty, then $D^b_{\operatorname{coh}}(Y_1)$ is the zero category, and so the equivalence forces $Y_2$ to be empty as well.
	Hence, we can assume each $Y_i$ is nonempty.
	In this case, each $Y_i$ is separated and of finite type over $k$, and hence, finitely presented over $k$.
	By \eqref{eq:depth_spectra} and \Cref{lem:max_closed_depth},
	\begin{displaymath}
		\dim Y_1
		=\max_{x\in|Y_1|\textrm{ closed}}\operatorname{depth}_x(Y_1)
		=\max_{y\in|Y_2|\textrm{ closed}}\operatorname{depth}_y(Y_2)
		=\dim Y_2.
	\end{displaymath}
	This completes the proof.
\end{proof}

\begin{definition}
	\label{def:orthogonal_decomposition}
	Let $\mathcal{T}$ be a triangulated category.
	An \textbf{orthogonal decomposition} of $\mathcal{T}$ is a finite collection $\{\mathcal{A}_i\}_{i\in I}$ of strictly full triangulated subcategories of $\mathcal{T}$ satisfying the following conditions:
	\begin{enumerate}
		\item If $i,i^\prime\in I$ are distinct, $A\in\mathcal{A}_i$, and $A^\prime\in\mathcal{A}_{i^\prime}$, then $\operatorname{Hom}_{\mathcal{T}}(A,A^\prime)=0$.
		\item If $E \in \mathcal{T}$, then there exists $A_i\in\mathcal{A}_i$ such that $E\cong \oplus_{i\in I}A_i$.
		\end{enumerate}
	We allow $I=\varnothing$, so that the zero category admits the empty orthogonal decomposition.
	In such a case, we denote this by $\mathcal{T}=\oplus_{i\in I}\mathcal{A}_i$.
	We say that $\mathcal{T}$ is \textbf{indecomposable} if it has a nonzero object and for every orthogonal decomposition $\mathcal{T}=\oplus_{i\in I}\mathcal{A}_i$, there exists only one $\mathcal{A}_i$ with a nonzero object.
\end{definition}

\begin{lemma}
	\label{lem:orthogonal_basic}
	Let $\mathcal{T}=\oplus_{i\in I}\mathcal{A}_i$ be an orthogonal decomposition of a triangulated category.
	\begin{enumerate}
		\item \label{lem:orthogonal_basic1}
		Each $\mathcal{A}_i$ is closed under direct summands in $\mathcal{T}$.
		\item \label{lem:orthogonal_basic2}
		If $i,i^\prime\in I$ are distinct, then $\mathcal{A}_i\cap\mathcal{A}_{i^\prime}$ consists only of zero objects.
		\item \label{lem:orthogonal_basic3}
		Let $F\colon\mathcal{T}\to\mathcal{T}^\prime$ be an exact equivalence.
		For each $i\in I$, denote by $F(\mathcal{A}_i)$ the strictly full subcategory of $\mathcal{T}^\prime$ consisting of objects isomorphic to $F(A)$ for some $A\in\mathcal{A}_i$.
		Then $\mathcal{T}^\prime=\oplus_{i\in I}F(\mathcal{A}_i)$ is an orthogonal decomposition.
		Moreover, $F(\mathcal{A}_i)$ is indecomposable if, and only if, $\mathcal{A}_i$ is indecomposable.
	\end{enumerate}
\end{lemma}

\begin{proof}
	We prove \eqref{lem:orthogonal_basic1}.
	Let $X\in\mathcal{A}_i$.
	Suppose that $X\cong Y\oplus Y^\prime$ in $\mathcal{T}$.
	Write $Y\cong\oplus_{l\in I}Y_l$ with $Y_l\in\mathcal{A}_l$.
	Fix $l\neq i$.
	Since $Y_l$ is a direct summand of $X$, the identity of $Y_l$ factors as a composite $Y_l\to X\to Y_l$.
	As $\operatorname{Hom}_{\mathcal{T}}(Y_l,X)=0$, the first morphism is zero.
	Hence, the identity of $Y_l$ is zero, and so $Y_l\cong0$.
	Thus, $Y\cong Y_i$.
	Since $\mathcal{A}_i$ is strictly full, $Y\in\mathcal{A}_i$.

	We prove \eqref{lem:orthogonal_basic2}.
	If $X\in\mathcal{A}_i\cap\mathcal{A}_{i^\prime}$, then $\operatorname{Hom}_{\mathcal{T}}(X,X)=0$.
	Hence, the identity of $X$ is zero, and so $X\cong0$.

	We prove \eqref{lem:orthogonal_basic3}.
	Since $F$ is an exact equivalence, each $F(\mathcal{A}_i)$ is a strictly full triangulated subcategory of $\mathcal{T}^\prime$.
	If $i\neq i^\prime$, $A\in\mathcal{A}_i$, and $A^\prime\in\mathcal{A}_{i^\prime}$, then
	\begin{displaymath}
		\operatorname{Hom}_{\mathcal{T}^\prime}(F(A),F(A^\prime))
		\cong\operatorname{Hom}_{\mathcal{T}}(A,A^\prime)
		=0.
	\end{displaymath}
	Note that every object of $\mathcal{T}^\prime$ is isomorphic to $F(X)$ for some $X\in\mathcal{T}$.
	Write $X\cong\oplus_{i\in I}A_i$ with $A_i\in\mathcal{A}_i$.
	As $F$ is additive, $F(X)\cong\oplus_{i\in I}F(A_i)$.
	Hence, $\mathcal{T}^\prime=\oplus_{i\in I}F(\mathcal{A}_i)$ is an orthogonal decomposition.

	Moreover, $F$ restricts to an exact equivalence from $\mathcal{A}_i$ to $F(\mathcal{A}_i)$.
	Applying the previous argument to this equivalence (with a quasi-inverse of it), it follows that: 1.\ orthogonal decompositions of $\mathcal{A}_i$ and of $F(\mathcal{A}_i)$ correspond to each other, and 2.\ an object is nonzero if, and only if, its image is nonzero.
	Consequently, $F(\mathcal{A}_i)$ is indecomposable if, and only if, $\mathcal{A}_i$ is indecomposable.
\end{proof}

\begin{lemma}
	\label{lem:orthogonal_unique}
	Let $\mathcal{T}$ be a triangulated category.
	Suppose that $\mathcal{T}=\oplus_{i\in I}\mathcal{A}_i$ and $\mathcal{T}=\oplus_{j\in J}\mathcal{B}_j$ are orthogonal decompositions.
	Assume that every $\mathcal{A}_i$ and every $\mathcal{B}_j$ is indecomposable.
	Then there exists a unique bijection $\sigma\colon I\to J$ where $\mathcal{A}_i=\mathcal{B}_{\sigma(i)}$ for every $i\in I$.
\end{lemma}

\begin{proof}
	Fix $i\in I$.
	For each $j\in J$, $\mathcal{A}_i\cap\mathcal{B}_j$ is a strictly full triangulated subcategory of $\mathcal{A}_i$.
	We claim that $\mathcal{A}_i=\oplus_{j\in J}(\mathcal{A}_i\cap\mathcal{B}_j)$ is an orthogonal decomposition.
	The required vanishing of morphisms holds.
	Indeed, they hold for $\mathcal{T}=\oplus_{j\in J}\mathcal{B}_j$.
	Pick $X\in\mathcal{A}_i$.
	Write $X\cong\oplus_{j\in J}X_j$ with $X_j\in\mathcal{B}_j$.
	Note that each $X_j$ is a direct summand of $X$.
	Hence, \Cref{lem:orthogonal_basic}\eqref{lem:orthogonal_basic1} implies $X_j\in\mathcal{A}_i$.
	Since $X_j\in\mathcal{A}_i\cap\mathcal{B}_j$, the claim follows.
	As $\mathcal{A}_i$ is indecomposable, there exists a unique $\sigma(i)\in J$ such that $\mathcal{A}_i\cap\mathcal{B}_{\sigma(i)}$ has a nonzero object.
	By the decomposition $X\cong\oplus_{j\in J}X_j$ above, it follows that $X_j\cong0$ for $j\neq\sigma(i)$.
	Hence, $X\cong X_{\sigma(i)}$, and so $X\in\mathcal{B}_{\sigma(i)}$.
	Indeed, $\mathcal{B}_{\sigma(i)}$ is strictly full.
	This shows that $\mathcal{A}_i\subseteq\mathcal{B}_{\sigma(i)}$.

	By symmetry, for every $j\in J$ there exists $\tau(j)\in I$ satisfying $\mathcal{B}_j\subseteq\mathcal{A}_{\tau(j)}$.
	Choose $i\in I$.
	Note that $\mathcal{A}_i\subseteq\mathcal{B}_{\sigma(i)}\subseteq\mathcal{A}_{\tau(\sigma(i))}$.
	Since $\mathcal{A}_i$ has a nonzero object, \Cref{lem:orthogonal_basic}\eqref{lem:orthogonal_basic2} implies $\tau(\sigma(i))=i$.
	Hence, $\mathcal{A}_i=\mathcal{B}_{\sigma(i)}$.
	Similarly, $\sigma(\tau(j))=j$ for every $j\in J$, and so $\sigma$ is a bijection.
	If $\sigma^\prime\colon I\to J$ is another such bijection, then $\mathcal{B}_{\sigma(i)}=\mathcal{A}_i=\mathcal{B}_{\sigma^\prime(i)}$ has a nonzero object.
	Hence, $\sigma(i)=\sigma^\prime(i)$ by \Cref{lem:orthogonal_basic}\eqref{lem:orthogonal_basic2}.
	This completes the proof.
\end{proof}

\begin{lemma}
	\label{lem:connected_components}
	Let $X$ be a Noetherian algebraic space.
	\begin{enumerate}
		\item \label{lem:connected_components1}
		$|X|$ has finitely many connected components, and each of them is open and closed in $|X|$.
		\item \label{lem:connected_components2}
		Let $j\colon X^\prime\to X$ be an open immersion such that $|j|(|X^\prime|)$ is closed in $|X|$.
		Then $X^\prime$ is Noetherian, the counit $\mathbf{L}j^\ast\mathbf{R}j_\ast A\to A$ is an isomorphism for every $A\in D_{\operatorname{qc}}(X^\prime)$, and $\mathbf{R}j_\ast$ restricts to a fully faithful exact functor $D^b_{\operatorname{coh}}(X^\prime)\to D^b_{\operatorname{coh}}(X)$.
		\item \label{lem:connected_components3}
		For each connected component of $|X|$, consider the essential image of the functor in \eqref{lem:connected_components2} associated to the corresponding open subspace of $X$.
		These essential images form an orthogonal decomposition of $D^b_{\operatorname{coh}}(X)$.
		\item \label{lem:connected_components4}
		If $|X|$ is nonempty and connected, then $D^b_{\operatorname{coh}}(X)$ is indecomposable.
	\end{enumerate}
\end{lemma}

\begin{proof}
	We prove \eqref{lem:connected_components1}.
	The Noetherian topological space $|X|$ has finitely many irreducible components.
	These irreducible components are connected.
	Hence, $|X|$ has finitely many connected components.
	Since connected components are closed, their complements are finite unions of closed subsets, so they are also open \cite[\href{https://stacks.math.columbia.edu/tag/004T}{Tag 004T}]{StacksProject}.

	We prove \eqref{lem:connected_components2}.
	By \cite[\href{https://stacks.math.columbia.edu/tag/03BZ}{Tag 03BZ}]{StacksProject}, the complement of $|j|(|X^\prime|)$ is an open subspace $j^{\prime\prime}\colon X^{\prime\prime}\to X$.
	Then $j$ and $j^{\prime\prime}$ imply that $X$ is the disjoint union of $X^\prime$ and $X^{\prime\prime}$.
	Since they are open and closed in $X$, both subspaces are Noetherian.
	Hence, restriction to these subspaces yields equivalences
	\begin{displaymath}
		\begin{aligned}
			D_{\operatorname{qc}}(X)&\cong D_{\operatorname{qc}}(X^\prime)\bigoplus D_{\operatorname{qc}}(X^{\prime\prime}),\\
			D^b_{\operatorname{coh}}(X)&\cong D^b_{\operatorname{coh}}(X^\prime)\bigoplus D^b_{\operatorname{coh}}(X^{\prime\prime}).
		\end{aligned}
	\end{displaymath}
	(see also \cite[Lemma A.4]{GuisadoVillaalgordo/Lank:2026} for the first equivalence).
	Under these equivalences, $\mathbf{L}j^\ast$ is projection onto the first factor and $\mathbf{R}j_\ast$ is its fully faithful inclusion.
	Consequently, the counit $\mathbf{L}j^\ast\mathbf{R}j_\ast A\to A$ is an isomorphism and $\mathbf{R}j_\ast$ preserves bounded pseudocoherent complexes.

	We prove \eqref{lem:connected_components3}.
	By \eqref{lem:connected_components1}, $X$ is the finite disjoint union of the open subspaces corresponding to its connected components.
	Moreover, the same decomposition of $D^b_{\operatorname{coh}}(X)$ as above yields the desired orthogonal decomposition.

	We prove \eqref{lem:connected_components4}.
	Let $D^b_{\operatorname{coh}}(X)=\oplus_{l\in L}\mathcal{A}_l$ be an orthogonal decomposition.
	Since $X$ is nonempty, $\mathcal{O}_X$ is nonzero.
	Its endomorphism ring is $\Gamma(X,\mathcal{O}_X)$.
	Note that the only idempotents are $0$ and $1$. 
	Indeed, an idempotent determines an open-and-closed partition of $|X|$, whereas $|X|$ is connected.
		Hence, $\mathcal{O}_X$ is indecomposable.
	Moreover, it lies in a unique $\mathcal{A}_{l_0}$: the decomposition writes $\mathcal{O}_X$ as a direct sum of objects of the $\mathcal{A}_l$, and a direct sum of two nonzero objects is not indecomposable, so exactly one summand is nonzero, and it contains $\mathcal{O}_X$ because $\mathcal{A}_{l_0}$ is strictly full.
	Pick a closed point $x\in|X|$.
	It exists because $|X|$ is nonempty and quasi-compact, and it is Kolmogorov as $X$ is decent \cite[\href{https://stacks.math.columbia.edu/tag/03JX}{Tags 03JX}, \href{https://stacks.math.columbia.edu/tag/03K6}{03K6}, and \href{https://stacks.math.columbia.edu/tag/005E}{005E}]{StacksProject}.
	Since $X$ is Noetherian, the residual space of $x$ is a closed immersion $i_x\colon\operatorname{Spec}(\kappa(x))\to X$ \cite[\href{https://stacks.math.columbia.edu/tag/0H1R}{Tags 0H1R} and \href{https://stacks.math.columbia.edu/tag/0H1U}{0H1U}]{StacksProject}.
	Moreover, the endomorphism ring of $\mathbf{R}(i_x)_\ast\mathcal{O}_{\operatorname{Spec}(\kappa(x))}$ is $\kappa(x)$ \cite[\href{https://stacks.math.columbia.edu/tag/08AM}{Tag 08AM}]{StacksProject}.
		Hence, the same argument shows that this object is indecomposable and lies in a unique $\mathcal{A}_l$.
	Hence, orthogonality forces $l=l_0$, since
	\begin{displaymath}
		\operatorname{Hom}_X(\mathcal{O}_X,\mathbf{R}(i_x)_\ast\mathcal{O}_{\operatorname{Spec}(\kappa(x))})\cong\kappa(x)\neq0.
	\end{displaymath}
	Now, suppose that $F\in\mathcal{A}_l$ is nonzero for some $l\neq l_0$.
	Set $b=\sup(F)$.
	Choose a point $x$ of the support of $\mathcal{H}^b(F)$ which is closed in this support.
	It exists because this support is nonempty and closed \cite[\href{https://stacks.math.columbia.edu/tag/07TZ}{Tag 07TZ}]{StacksProject}, and hence, quasi-compact and Kolmogorov \cite[\href{https://stacks.math.columbia.edu/tag/005E}{Tag 005E}]{StacksProject}.
	As the support is closed in $|X|$, so is $x$.
	Note that $\mathbf{R}(i_x)_\ast\mathcal{O}_{\operatorname{Spec}(\kappa(x))}$ is a coherent sheaf placed in degree zero.
	Indeed, the pushforward along the closed immersion $i_x$ is exact \cite[\href{https://stacks.math.columbia.edu/tag/0A4K}{Tag 0A4K}]{StacksProject}.
	In particular, we can identify $\mathbf{R}(i_x)_\ast\mathcal{O}_{\operatorname{Spec}(\kappa(x))}$ with $(i_x)_\ast\mathcal{O}_{\operatorname{Spec}(\kappa(x))}$.
	Consider the standard truncation triangle
	\begin{displaymath}
		\tau^{\leq b-1}F\to F\to\mathcal{H}^b(F)[-b]\to(\tau^{\leq b-1}F)[1].
	\end{displaymath}
	As morphisms from objects of $D^{\leq b-1}(\mathcal{O}_X)$ to objects of $D^{\geq b}(\mathcal{O}_X)$ vanish, applying the functor
	\begin{displaymath}
		\operatorname{Hom}_X(-,\mathbf{R}(i_x)_\ast\mathcal{O}_{\operatorname{Spec}(\kappa(x))}[-b])
	\end{displaymath}
	to this triangle gives the first isomorphism below.
	The second holds because morphisms in $D(\mathcal{O}_X)$ between sheaves placed in degree zero are morphisms of $\mathcal{O}_X$-modules, combined with the adjunction between $i_x^\ast$ and $(i_x)_\ast$ on modules:
	\begin{displaymath}
		\begin{aligned}
			&\operatorname{Hom}_X(F,\mathbf{R}(i_x)_\ast\mathcal{O}_{\operatorname{Spec}(\kappa(x))}[-b])\\
			&\qquad\cong\operatorname{Hom}_X(\mathcal{H}^b(F),\mathbf{R}(i_x)_\ast\mathcal{O}_{\operatorname{Spec}(\kappa(x))})\\
			&\qquad\cong\operatorname{Hom}(i_x^\ast\mathcal{H}^b(F),\mathcal{O}_{\operatorname{Spec}(\kappa(x))})\neq0.
		\end{aligned}
	\end{displaymath}
	Note that the last nonvanishing follows because $i_x^\ast\mathcal{H}^b(F)$ is a nonzero $\kappa(x)$-vector space \cite[\href{https://stacks.math.columbia.edu/tag/07TZ}{Tag 07TZ}]{StacksProject}.
	However, this contradicts orthogonality because the second argument belongs to $\mathcal{A}_{l_0}$.
	Therefore, $\mathcal{A}_l$ is zero for every $l\neq l_0$, which completes the proof.

	Statement \eqref{lem:connected_components4} is the analogue for Noetherian algebraic spaces of \cite[Example 3.2]{Bridgeland:1999}, which proves the corresponding assertion for the bounded derived category of a scheme.
\end{proof}

\begin{lemma}
	\label{lem:transform_restriction}
	Let $f_i\colon Y_i\to\operatorname{Spec}(k)$ be proper morphisms of algebraic spaces with $k$ a field.
	Suppose that $j_C\colon C\to Y_1$ and $j_D\colon D\to Y_2$ are open immersions such that $|j_C|(|C|)$ and $|j_D|(|D|)$ are closed.
	Then $C$ and $D$ are proper over $k$.
	Moreover, for every $E\in D_{\operatorname{qc}}(Y_1\times_kY_2)$ and $A\in D^b_{\operatorname{coh}}(C)$, there exists a natural isomorphism
	\begin{displaymath}
		\mathbf{L}(j_D)^\ast\Phi_E(\mathbf{R}(j_C)_\ast A)
		\cong\Phi_{\mathbf{L}(j_C\times_k j_D)^\ast E}(A),
	\end{displaymath}
	where the right-hand side is the integral transform with kernel $\mathbf{L}(j_C\times_k j_D)^\ast E\in D_{\operatorname{qc}}(C\times_kD)$.
\end{lemma}

\begin{proof}
	By the hypotheses, we can identify $C$ and $D$ with open-and-closed subspaces of $Y_1$ and $Y_2$.
	Note that each inclusion is also a closed immersion.
	Indeed, this may be checked locally on the target.
	Hence, $j_C$ and $j_D$ are proper, and so are $f_1\circ j_C$ and $f_2\circ j_D$.

		We prove the desired natural isomorphism.
	Write $p_i\colon Y_1\times_kY_2\to Y_i$ for the projections, $p_1^\prime\colon Y_1\times_kD\to Y_1$ and $p_2^\prime\colon Y_1\times_kD\to D$ for the intermediate projections, and $q_1\colon C\times_kD\to C$ and $q_2\colon C\times_kD\to D$ for the remaining projections.
	Consider the fibered squares
	\begin{displaymath}
		\begin{tikzcd}
			{Y_1\times_kD} & {Y_1\times_kY_2} && {C\times_kD} & {Y_1\times_kD} \\
			D & {Y_2} && C & {Y_1.}
			\arrow["{\operatorname{id}_{Y_1}\times_kj_D}", from=1-1, to=1-2]
			\arrow["{p_2^\prime}"', from=1-1, to=2-1]
			\arrow["{p_2}", from=1-2, to=2-2]
			\arrow["{j_D}"', from=2-1, to=2-2]
			\arrow["{j_C\times_k\operatorname{id}_D}", from=1-4, to=1-5]
			\arrow["{q_1}"', from=1-4, to=2-4]
			\arrow["{p_1^\prime}", from=1-5, to=2-5]
			\arrow["{j_C}"', from=2-4, to=2-5]
		\end{tikzcd}
	\end{displaymath}
	Both are Tor-independent, the first because $j_D$ is flat and the second because $p_1^\prime$ is flat; the remaining hypotheses of \cite[\href{https://stacks.math.columbia.edu/tag/08IR}{Tag 08IR}]{StacksProject} are quasi-compactness and quasi-separatedness of the morphisms involved, which holds since the spaces are Noetherian.
	By \Cref{lem:connected_components}\eqref{lem:connected_components2}, $\mathbf{R}(j_C)_\ast A\in D^b_{\operatorname{coh}}(Y_1)$.
	Since $j_D$ is flat and $p_2$ is proper, \cite[\href{https://stacks.math.columbia.edu/tag/08IR}{Tag 08IR}]{StacksProject} and monoidality of derived pullback yields
	\begin{displaymath}
		\begin{aligned}
			\mathbf{L}(j_D)^\ast\Phi_E(\mathbf{R}(j_C)_\ast A)
			&\cong\mathbf{R}(p_2^\prime)_\ast\mathbf{L}(\operatorname{id}_{Y_1}\times_k j_D)^\ast
			(\mathbf{L}(p_1)^\ast\mathbf{R}(j_C)_\ast A\otimes^{\mathbf{L}}E)\\
			&\cong\mathbf{R}(p_2^\prime)_\ast
			(\mathbf{L}(p_1^\prime)^\ast\mathbf{R}(j_C)_\ast A
			\otimes^{\mathbf{L}}\mathbf{L}(\operatorname{id}_{Y_1}\times_k j_D)^\ast E).
		\end{aligned}
	\end{displaymath}
	Here we used $p_1\circ(\operatorname{id}_{Y_1}\times_k j_D)=p_1^\prime$.
	Note that $p_1^\prime$ is flat
	and $j_C$ is quasi-compact quasi-separated.
	Applying flat base change to the fibered square defined by $j_C$ and $p_1^\prime$ gives
	\begin{displaymath}
		\mathbf{L}(p_1^\prime)^\ast\mathbf{R}(j_C)_\ast A
		\cong\mathbf{R}(j_C\times_k\operatorname{id}_D)_\ast\mathbf{L}(q_1)^\ast A.
	\end{displaymath}
	Note that $j_C\times_k\operatorname{id}_D$ is quasi-compact quasi-separated.
	By \cite[\href{https://stacks.math.columbia.edu/tag/08IN}{Tag 08IN}]{StacksProject}, 
	\begin{displaymath}
		\begin{aligned}
			&\mathbf{R}(j_C\times_k\operatorname{id}_D)_\ast\mathbf{L}(q_1)^\ast A
			\otimes^{\mathbf{L}}\mathbf{L}(\operatorname{id}_{Y_1}\times_k j_D)^\ast E\\
			&\qquad\cong\mathbf{R}(j_C\times_k\operatorname{id}_D)_\ast
			(\mathbf{L}(q_1)^\ast A\otimes^{\mathbf{L}}
			\mathbf{L}(j_C\times_k\operatorname{id}_D)^\ast
			\mathbf{L}(\operatorname{id}_{Y_1}\times_k j_D)^\ast E)\\
			&\qquad\cong\mathbf{R}(j_C\times_k\operatorname{id}_D)_\ast
			(\mathbf{L}(q_1)^\ast A\otimes^{\mathbf{L}}\mathbf{L}(j_C\times_k j_D)^\ast E).
		\end{aligned}
	\end{displaymath}
	Here the last isomorphism uses $(\operatorname{id}_{Y_1}\times_k j_D)\circ(j_C\times_k\operatorname{id}_D)=j_C\times_k j_D$.
	Lastly, $p_2^\prime\circ(j_C\times_k\operatorname{id}_D)=q_2$, and so \cite[\href{https://stacks.math.columbia.edu/tag/0D6E}{Tag 0D6E}]{StacksProject} gives
	\begin{displaymath}
		\begin{aligned}
			\mathbf{L}(j_D)^\ast\Phi_E(\mathbf{R}(j_C)_\ast A)
			&\cong\mathbf{R}(q_2)_\ast
			(\mathbf{L}(q_1)^\ast A\otimes^{\mathbf{L}}\mathbf{L}(j_C\times_k j_D)^\ast E)\\
			&=\Phi_{\mathbf{L}(j_C\times_k j_D)^\ast E}(A).
		\end{aligned}
	\end{displaymath}
	As all the isomorphisms are natural in $A$, we complete the proof.
\end{proof}

\begin{corollary}
	\label{cor:components_categorical}
	Let $X_1$ and $X_2$ be Noetherian algebraic spaces and let $F\colon D^b_{\operatorname{coh}}(X_1)\to D^b_{\operatorname{coh}}(X_2)$ be an exact equivalence.
	Then there exists a unique bijection between the connected components of $|X_1|$ and those of $|X_2|$ such that, for matched components $X_1^\prime$ and $X_2^\prime$, the equivalence $F$ restricts, under the identifications of \Cref{lem:connected_components}\eqref{lem:connected_components3}, to an equivalence $D^b_{\operatorname{coh}}(X_1^\prime)\to D^b_{\operatorname{coh}}(X_2^\prime)$.
\end{corollary}

\begin{proof}
	By \Cref{lem:connected_components}\eqref{lem:connected_components1}, \eqref{lem:connected_components3}, and \eqref{lem:connected_components4}, each $D^b_{\operatorname{coh}}(X_i)$ decomposes orthogonally into indecomposable pieces indexed by the connected components of $|X_i|$.
	By \Cref{lem:orthogonal_basic}\eqref{lem:orthogonal_basic3} and \Cref{lem:orthogonal_unique}, the equivalence $F$ determines a unique bijection between the two sets of components and identifies the corresponding pieces.
	Composing with the equivalences of \Cref{lem:connected_components}\eqref{lem:connected_components3} gives the assertion.
\end{proof}

\begin{proposition}
	\label{prop:connected_components}
	Let $f_i\colon Y_i\to\operatorname{Spec}(k)$ be proper morphisms of algebraic spaces with $k$ a field.
	Assume that $Y_1$ and $Y_2$ are Fourier--Mukai partners given by $E\in D^b_{\operatorname{coh}}(Y_1\times_kY_2)$.
	Then there exists a bijection between the connected components of $|Y_1|$ and those of $|Y_2|$ with the following property: If $j_C\colon C\to Y_1$ and $j_D\colon D\to Y_2$ are the open immersions of the open subspaces corresponding to matched connected components, then $C$ and $D$ are proper over $k$, they are Fourier--Mukai partners given by $\mathbf{L}(j_C\times_k j_D)^\ast E$, and
	\begin{displaymath}
		\dim C=\dim D,
		\qquad
		\{\operatorname{depth}_c(C):c\in|C|\textrm{ closed}\}
		=\{\operatorname{depth}_d(D):d\in|D|\textrm{ closed}\}.
	\end{displaymath}
\end{proposition}

\begin{proof}
		Each $Y_i$ is Noetherian because it is proper over the field $k$.
	If one of the $Y_i$ is empty, then so is the other, as an equivalence forces both derived categories to vanish, and there is no matched pair to consider.
	By \Cref{cor:components_categorical}, applied to the restriction of $\Phi_E$ to $D^b_{\operatorname{coh}}$, there exists a bijection between the connected components of $|Y_1|$ and those of $|Y_2|$ such that $\Phi_E$ maps the summand of the orthogonal decomposition of \Cref{lem:connected_components}\eqref{lem:connected_components3} corresponding to a component of $|Y_1|$ onto the summand corresponding to the matched component of $|Y_2|$.

	Fix matched components with open immersions $j_C\colon C\to Y_1$ and $j_D\colon D\to Y_2$.
	By \Cref{lem:connected_components}\eqref{lem:connected_components2}, the composite
	\begin{displaymath}
		D^b_{\operatorname{coh}}(C)
		\xrightarrow{\mathbf{R}(j_C)_\ast}D^b_{\operatorname{coh}}(Y_1)
		\xrightarrow{\Phi_E}D^b_{\operatorname{coh}}(Y_2)
		\xrightarrow{\mathbf{L}(j_D)^\ast}D^b_{\operatorname{coh}}(D)
	\end{displaymath}
		is an equivalence.
	Indeed, by the choice of the bijection, $\Phi_E$ maps the essential image of $\mathbf{R}(j_C)_{\ast}$ onto that of $\mathbf{R}(j_D)_{\ast}$.
	Moreover, by \Cref{lem:connected_components}\eqref{lem:connected_components2}, $\mathbf{R}(j_C)_{\ast}$ is an equivalence onto the former, and $\mathbf{L}(j_D)^{\ast}$ restricts to an equivalence from the latter onto $D^b_{\operatorname{coh}}(D)$.
	By \Cref{lem:transform_restriction}, $C$ and $D$ are proper over $k$, and this composite is naturally isomorphic to $\Phi_{\mathbf{L}(j_C\times_k j_D)^\ast E}$.
	Since $j_C\times_k j_D$ is flat, the associated kernel belongs to $D^b_{\operatorname{coh}}(C\times_kD)$.
	Consequently, $C$ and $D$ are Fourier--Mukai partners, \Cref{cor:dimension_invariance} gives $\dim C=\dim D$, and \Cref{thm:depth_correspondence} applied to this kernel gives the asserted equality of the sets of depths at closed points.
\end{proof}

\begin{lemma}
	\label{lem:equidimensional_closed_points}
	Let $X$ be a nonempty algebraic space finitely presented over a field $k$.
	If there exists an integer $e$ such that $\dim_z^{\mathrm{loc}}(X)=e$ for every closed point $z\in|X|$, then $\dim X=e$ and $X$ is equidimensional.
\end{lemma}

\begin{proof}
	Fix $x\in|X|$.
	Choose an \'{e}tale germ $(U,u)\to(X,x)$ with $U$ affine.
	For each irreducible component $Z$ of $U$ containing $u$, \Cref{lem:dimension_charts}\eqref{lem:dimension_charts3} yields some $v\in Z$ mapping to a closed point $z\in|X|$ where $\dim(\mathcal{O}_{U,v})=\dim Z$.
	By \Cref{def:local-depth-dimension}, the hypothesis implies $\dim Z=\dim_z^{\mathrm{loc}}(X)=e$.
	Hence, every irreducible component of $U$ containing $u$ has dimension $e$.
	Applying \cite[\href{https://stacks.math.columbia.edu/tag/0A21}{Tag 0A21}]{StacksProject}, it follows that $\dim_x(X)=\dim_u(U)=e$.
	As $x$ was arbitrary, $\dim X=e$ and $X$ is equidimensional.
\end{proof}

\begin{lemma}
	\label{lem:connected_cm_equidimensional}
	Let $X$ be a nonempty algebraic space finitely presented over a field $k$.
	If $|X|$ is connected and $X$ is Cohen--Macaulay, then $X$ is equidimensional.
\end{lemma}

\begin{proof}
	Pick a closed point $x\in|X|$.
	Choose an \'{e}tale germ $a\colon(U,u)\to(X,x)$ with $U$ affine.
	By \Cref{lem:closed_points_charts}, $u$ is closed in $U$.
	Since $X$ is Cohen--Macaulay, \Cref{lem:local-depth-dimension} implies $\mathcal{O}_{U,u}$ is Cohen--Macaulay.
	By \cite[\href{https://stacks.math.columbia.edu/tag/0BUS}{Tag 0BUS}]{StacksProject}, every minimal prime $\mathfrak{p}$ of this local ring (which is an associated prime) satisfies $\dim(\mathcal{O}_{U,u}/\mathfrak{p})=\dim(\mathcal{O}_{U,u})$.
	As $u$ is closed in the finite type $k$-scheme $U$, each irreducible component $Z$ of $U$ containing $u$ has dimension
	\begin{displaymath}
		\dim Z=\dim(\mathcal{O}_{U,u})=\dim_x^{\mathrm{loc}}(X).
	\end{displaymath}
	We can remove from $U$ the irreducible components not containing $u$.
	Denote the resulting open neighborhood of $u$ by $V$.
	If $z\in|a|(V)$ is closed in $|X|$, pick some $v\in V$ mapping to $z$.
		Since $v\in V$, and $V$ was obtained from $U$ by removing the irreducible components of $U$ not containing $u$, every irreducible component of $U$ containing $v$ also contains $u$.
	By \Cref{lem:closed_points_charts}, $v$ is closed in $U$.
	Consequently, \cite[\href{https://stacks.math.columbia.edu/tag/0A21}{Tag 0A21}]{StacksProject} ensures that
	\begin{displaymath}
		\dim_z^{\mathrm{loc}}(X)=\dim(\mathcal{O}_{U,v})=\dim_x^{\mathrm{loc}}(X).
	\end{displaymath}
	Thus, each closed point $x$ has an open neighborhood on whose closed points the local dimension is $\dim_x^{\mathrm{loc}}(X)$.

	Note that these neighborhoods cover $|X|$.
	Indeed, every point specializes to a closed point, whereas any open neighborhood of that closed point contains its generalizations.
	Moreover, two such neighborhoods with nonempty intersection have the same constant value.
	Indeed, their intersection is a nonempty open subset of $|X|$, and so it contains a closed point of $|X|$, as $|X|$ is Jacobson \cite[\href{https://stacks.math.columbia.edu/tag/0BA3}{Tags 0BA3} and \href{https://stacks.math.columbia.edu/tag/0BA6}{0BA6}]{StacksProject}.
	Consequently, for each integer, the union of the neighborhoods with that constant value is open, these unions are pairwise disjoint, and they cover $|X|$.
	Since $|X|$ is connected, only one of them is nonempty.
	Consequently, $\dim_z^{\mathrm{loc}}(X)$ is independent of the closed point $z\in|X|$, and \Cref{lem:equidimensional_closed_points} finishes the proof.
\end{proof}

\begin{proposition}
	\label{prop:cm_invariance}
	Let $f_i\colon Y_i\to\operatorname{Spec}(k)$ be proper morphisms of algebraic spaces with $k$ a field.
	Assume that $Y_1$ and $Y_2$ are Fourier--Mukai partners given by some $E\in D^b_{\operatorname{coh}}(Y_1\times_k Y_2)$.
	Then $Y_1$ is Cohen--Macaulay if, and only if, $Y_2$ is Cohen--Macaulay.
	Moreover, $Y_1$ is Cohen--Macaulay and equidimensional if, and only if, $Y_2$ is such.
\end{proposition}

\begin{proof}
		The case of an empty partner is immediate: the equivalence forces both derived categories to be zero, hence both spaces to be empty, and then the assertions are vacuous.
	Hence, we may assume that $Y_1$ and $Y_2$ are nonempty.
	By \Cref{cor:dimension_invariance}, they have the same dimension $d$.
	We first prove the second claim. 
	To do so, we show that a nonempty algebraic space $X$ finitely presented over $k$ and of dimension $d$ is Cohen--Macaulay and equidimensional if, and only if, $\operatorname{depth}_z(X)=d$ for every closed point $z\in|X|$.
	Indeed, \Cref{lem:dimension_charts}\eqref{lem:dimension_charts2} gives
	\begin{displaymath}
		\operatorname{depth}_z(X)\leq\dim_z^{\mathrm{loc}}(X)=\dim_z(X)\leq d.
	\end{displaymath}
	If $X$ is Cohen--Macaulay and equidimensional, then these inequalities are equalities.
	Conversely, by \Cref{lem:cm_closed_points}, equality of the depth with $d$ at every closed point implies $X$ must be Cohen--Macaulay.
	Hence, $\dim_z^{\mathrm{loc}}(X)=d$ at every closed point.
	By \Cref{lem:equidimensional_closed_points}, $X$ is equidimensional.
		Since $Y_1$ and $Y_2$ have the same dimension and, by \eqref{eq:depth_spectra}, the same set of depths at closed points, the characterization above applies to both or to neither; this proves the second claim.

	For the first claim, let $Y_1$ be Cohen--Macaulay.
		Note that every connected component $C$ of $Y_1$ is Cohen--Macaulay and equidimensional.
	Indeed, Cohen--Macaulayness is local, so $C$ is Cohen--Macaulay, and $C$ is equidimensional by \Cref{lem:connected_cm_equidimensional}.
	Applying \Cref{prop:connected_components}, the connected components of $Y_1$ and $Y_2$ are matched in Fourier--Mukai pairs.
	By the second claim, every connected component of $Y_2$ is Cohen--Macaulay.
	As these components are open and cover $Y_2$, the space $Y_2$ is Cohen--Macaulay.
	The converse follows by symmetry.
\end{proof}

\begin{proof}
	[Proof of \Cref{cor:intro_relative}]
	Let $E\in D^b_{\operatorname{coh}}(Y_1\times_S Y_2)$ be a kernel giving the partnership.
	By \cite[Remark 2.7]{GuisadoVillaalgordo/Lank:2026}, $E$ is relatively perfect over both factors.
	We first record a consequence of base change.
	Let $t\colon T\to S$ be a morphism where $T$ is the spectrum of a field.
	Since $S$ is Noetherian, it is quasi-separated.
	By \cite[\href{https://stacks.math.columbia.edu/tag/09TF}{Tag 09TF}]{StacksProject}, $t$ is affine, and hence, quasi-affine.
	Denote by $t^\prime\colon Y_1\times_S Y_2\times_S T\to Y_1\times_S Y_2$ the base change of $t$.
	Note that $Y_1\times_S Y_2\times_S T$ is isomorphic to $(Y_1\times_S T)\times_T(Y_2\times_S T)$.
	By \cite[Corollary 6.18]{GuisadoVillaalgordo/Lank:2026}, $\mathbf{L}(t^\prime)^\ast E$ is relatively perfect over both factors.
	In particular, it has bounded cohomology \cite[Remark 4.2(1)]{GuisadoVillaalgordo/Lank:2026}.
	It is also pseudocoherent \cite[\href{https://stacks.math.columbia.edu/tag/08H4}{Tag 08H4}]{StacksProject}.
	Hence, $\mathbf{L}(t^\prime)^\ast E$ belongs to $D^b_{\operatorname{coh}}((Y_1\times_S T)\times_T(Y_2\times_S T))$.
	By \cite[Theorem 10.13]{GuisadoVillaalgordo/Lank:2026}, the integral transform with kernel $\mathbf{L}(t^\prime)^\ast E$ is an equivalence on $D^b_{\operatorname{coh}}$.
	As $Y_1\times_S T$ and $Y_2\times_S T$ are proper and flat over $T$, being base changes of $f_1$ and $f_2$, it follows that
	\begin{equation}
		\label{eq:intro_relative_fibers}
		Y_1\times_S T\textrm{ and }Y_2\times_S T\textrm{ are Fourier--Mukai partners over }T.
	\end{equation}

	We prove \eqref{cor:intro_relative1} and \eqref{cor:intro_relative2}.
	As $S$ is quasi-separated, it is decent, and so every $s\in|S|$ is represented by a monomorphism $\operatorname{Spec}(\kappa(s))\to S$ \cite[\href{https://stacks.math.columbia.edu/tag/03JX}{Tags 03JX}, \href{https://stacks.math.columbia.edu/tag/03K4}{03K4}, and \href{https://stacks.math.columbia.edu/tag/0EMW}{0EMW}]{StacksProject}.
	Applying \eqref{eq:intro_relative_fibers} to this morphism, \Cref{introthm:main} applies to $(Y_1)_s$ and $(Y_2)_s$ over $\kappa(s)$.

	We prove \eqref{cor:intro_relative3}.
	Let $y\in|Y_i|$ with image $s\in|S|$.
	By \cite[\href{https://stacks.math.columbia.edu/tag/03H4}{Tag 03H4}]{StacksProject}, there exists $y^\prime\in|(Y_i)_s|$ mapping to $y$, and every point of $|(Y_i)_s|$ maps to a point of $|Y_i|$ lying over $s$.
	By \cite[\href{https://stacks.math.columbia.edu/tag/04NS}{Tag 04NS}]{StacksProject}, the relative dimension of $f_i$ at $y$ equals that of $(Y_i)_s\to\operatorname{Spec}(\kappa(s))$ at $y^\prime$, which is $\dim_{y^\prime}((Y_i)_s)$ \cite[\href{https://stacks.math.columbia.edu/tag/04NH}{Tag 04NH}]{StacksProject}.
	Hence, by \cite[\href{https://stacks.math.columbia.edu/tag/06LR}{Tags 06LR} and \href{https://stacks.math.columbia.edu/tag/04N6}{04N6}]{StacksProject}, $f_i$ has relative dimension $\leq d$ if, and only if, $\dim (Y_i)_s\leq d$ for every $s\in|S|$.
	By \eqref{cor:intro_relative1}, the latter condition does not depend on $i$.

	Lastly, we prove \eqref{cor:intro_relative4} and \eqref{cor:intro_relative5}.
	Choose a surjective \'{e}tale morphism $V\to S$ where $V$ is a scheme; the same $V$ is used for both $i$.
	For each $i$, choose a surjective \'{e}tale morphism $U_i\to Y_i\times_S V$ where $U_i$ is a scheme, and let $h_i\colon U_i\to V$ be its composite with the projection.
	Note that $h_i$ is flat, being the composite of an \'{e}tale morphism with a base change of $f_i$.
	Moreover, the composite $U_i\to Y_i$ is surjective and \'{e}tale, as $Y_i\times_S V\to Y_i$ is a base change of $V\to S$.
	As $V$ and $U_i$ are locally Noetherian, all local rings below are Noetherian.
	Fix $v\in V$.
	Denote by $\mathfrak{m}_v$ the maximal ideal of $\mathcal{O}_{V,v}$.
	For every $u\in U_i$ with $h_i(u)=v$, the local homomorphism $\mathcal{O}_{V,v}\to\mathcal{O}_{U_i,u}$ is flat, and $\mathcal{O}_{U_i,u}/\mathfrak{m}_v\mathcal{O}_{U_i,u}$ is the local ring of $U_i\times_V\operatorname{Spec}(\kappa(v))$ at $u$.
	By \cite[\href{https://stacks.math.columbia.edu/tag/00ON}{Tags 00ON} and \href{https://stacks.math.columbia.edu/tag/045J}{045J}]{StacksProject},
	\begin{displaymath}
		\dim(\mathcal{O}_{U_i,u})=\dim(\mathcal{O}_{V,v})+\dim(\mathcal{O}_{U_i\times_V\operatorname{Spec}(\kappa(v)),u}),
	\end{displaymath}
	and $\mathcal{O}_{U_i,u}$ is Cohen--Macaulay if, and only if, both $\mathcal{O}_{V,v}$ and $\mathcal{O}_{U_i\times_V\operatorname{Spec}(\kappa(v)),u}$ are Cohen--Macaulay.
	Moreover, the base change $U_i\times_V\operatorname{Spec}(\kappa(v))\to Y_i\times_S\operatorname{Spec}(\kappa(v))$ of $U_i\to Y_i\times_S V$ is surjective and \'{e}tale.
	Note that the points of $U_i\times_V\operatorname{Spec}(\kappa(v))$ are the points $u\in U_i$ with $h_i(u)=v$.
	Hence, by \Cref{def:cohen_macaulay}, \Cref{lem:local-depth-dimension}, and \cite[\href{https://stacks.math.columbia.edu/tag/04N5}{Tags 04N5}, \href{https://stacks.math.columbia.edu/tag/04N6}{04N6}, and \href{https://stacks.math.columbia.edu/tag/04MU}{04MU}]{StacksProject}:
	\begin{itemize}
		\item $Y_i\times_S\operatorname{Spec}(\kappa(v))$ is empty if, and only if, no point of $U_i$ maps to $v$;
		\item the Krull dimension of $Y_i\times_S\operatorname{Spec}(\kappa(v))$ is the supremum of $\dim(\mathcal{O}_{U_i\times_V\operatorname{Spec}(\kappa(v)),u})$ over all $u\in U_i$ mapping to $v$;
		\item $Y_i\times_S\operatorname{Spec}(\kappa(v))$ is Cohen--Macaulay if, and only if, each of these local rings is Cohen--Macaulay.
	\end{itemize}
	By the same reasoning applied to $U_i\to Y_i$, $\dim Y_i$ is the supremum of $\dim(\mathcal{O}_{U_i,u})$ over $u\in U_i$, and $Y_i$ is Cohen--Macaulay if, and only if, each $\mathcal{O}_{U_i,u}$ is Cohen--Macaulay.
	Combining yields
	\begin{displaymath}
		\dim Y_i=\sup\{\dim(\mathcal{O}_{V,v})+\dim(Y_i\times_S\operatorname{Spec}(\kappa(v))) : v\in V\textrm{ and }Y_i\times_S\operatorname{Spec}(\kappa(v))\neq\varnothing\}.
	\end{displaymath}
	Similarly, $Y_i$ is Cohen--Macaulay if, and only if, $\mathcal{O}_{V,v}$ and $Y_i\times_S\operatorname{Spec}(\kappa(v))$ are Cohen--Macaulay for every $v\in V$ such that $Y_i\times_S\operatorname{Spec}(\kappa(v))$ is nonempty.

	Fix $v\in V$.
	Apply \eqref{eq:intro_relative_fibers} to the composite $\operatorname{Spec}(\kappa(v))\to V\to S$.
	An algebraic space is empty if, and only if, its bounded derived category is zero.
	Hence, $Y_1\times_S\operatorname{Spec}(\kappa(v))$ is empty if, and only if, $Y_2\times_S\operatorname{Spec}(\kappa(v))$ is empty.
	If they are nonempty, \Cref{cor:dimension_invariance,prop:cm_invariance} show that they have the same Krull dimension, and that one is Cohen--Macaulay if, and only if, the other is.
	As $V$ does not depend on $i$, the two descriptions above do not depend on $i$.
	This completes the proof.
\end{proof}

\bibliographystyle{alpha}
\bibliography{mainbib}

\end{document}